\documentclass[11pt]{article}

\usepackage{amssymb,latexsym,amsfonts,verbatim,amscd}

\newtheorem{Def}{Definition}[section]
\newtheorem{Thm}[Def]{Theorem}
\newtheorem{Lem}[Def]{Lemma}
\newtheorem{Prop}[Def]{Proposition}
\newtheorem{Cor}[Def]{Corollary}
\newtheorem{Rem}[Def]{Remark}

\newtheorem{Quest}[Def]{Question}
\font\nat msbm10 scaled\magstephalf
\def\N{\hbox{\nat\char78}}

\def\P{\hbox{\Nat\char80}}

\font\mata=msam10 
\def\restr{\mbox{\mata\char22}}
\def\telos{\hfill$\dashv$}
  
\font\Nat=msbm10

\begin{document}
\sloppy

\title{Large transitive models  in local {\rm ZFC}}
\author{Athanassios Tzouvaras}

\date{}
\maketitle

\begin{center}
 Department  of Mathematics\\  Aristotle University of Thessaloniki \\
 541 24 Thessaloniki, Greece \\
  e-mail: \verb"tzouvara@math.auth.gr"
\end{center}

\abstract{This paper is a sequel to \cite{Tz10},  where a local version of ZFC, LZFC, was introduced and examined and transitive models of ZFC with properties that resemble  large cardinal properties, namely Mahlo and $\Pi_1^1$-indescribable models, were considered. By analogy we refer to such models as ``large models'', and the properties in question as ``large model properties''. Continuing here in the same spirit  we consider further large model properties, that resemble  stronger large cardinals, namely, ``elementarily embeddable'',  ``extendible'' and ``strongly extendible'', ``critical'' and ``strongly critical'', ``self-critical'' and ``strongly self-critical'', the definitions of which involve   elementary embeddings. Each large model property $\phi$ gives rise to a localization  axiom $Loc^{\phi}({\rm ZFC})$ saying that every set belongs to a transitive model of ZFC satisfying $\phi$. The theories ${\rm LZFC}^\phi={\rm LZFC}$+$Loc^{\phi}({\rm ZFC})$ are local analogues of the theories ZFC+``there is a proper class of large cardinals $\psi$'', where $\psi$ is a large cardinal property. If $sext(x)$ is the  property of  strong extendibility, it is shown that  ${\rm LZFC}^{sext}$ proves Powerset and $\Sigma_1$-Collection. In order to refute $V=L$ over LZFC, we  combine the existence of strongly critical models with an axiom of different flavor,  the Tall Model Axiom ($TMA$). $V=L$ can also be refuted by  $TMA$ plus the axiom $GC$ saying that ``there is a greatest cardinal'', although it is not known if $TMA+GC$ is consistent over LZFC. Finally Vop\v{e}nka's Principle ($VP$)  and its impact on LZFC are examined. It is shown that   ${\rm LZFC}^{sext}+VP$ proves Powerset and Replacement, i.e., ZFC is fully recovered. The same is true for some weaker variants of ${\rm LZFC}^{sext}$. Moreover the theories LZFC$^{sext}$+$VP$ and ZFC+$VP$ are shown to be identical.}

\vskip 0.2in


\vskip 0.2in


\section{Introduction}
In \cite{Tz10} we initiated the study of a local version of ZFC called LZFC, whose main  axiom, denoted $Loc({\rm ZFC})$, says that ``every set belongs to a  transitive model of ZFC''. Namely, this is the statement:
\begin{equation} \label{E:fund}
(Loc({\rm ZFC}))
\hspace{.5\columnwidth minus .5\columnwidth} \forall x\exists
y(x\in y \ \wedge \ Tr(y) \ \wedge \ (y,\in)\models{\rm
ZFC}).\hspace{.5\columnwidth minus .5\columnwidth} \llap{}
\end{equation}
LZFC is the theory BST+$Loc({\rm ZFC})$, where BST (which stands for ``basic set theory'') is the set of the following elementary axioms: Empty set, Extensionality, Pair, Union, Cartesian Product, ``$\omega$ exists'', $\Delta_0$-Separation. Actually BST is needed only for the formalization of the notion ``model of ZFC'' involved in $Loc({\rm ZFC})$ (see Remark 2.6 of \cite{Tz10}).

LZFC lacks the Powerset and Replacement axioms, as well as $\in$-induction.\footnote{The usual Regularity axiom holds in LZFC. $\in$-induction and, equivalently, $On$-induction is missing because full Separation is missing.}  These principles  hold only locally in transitive  set models which, in compensation, exist everywhere across the universe. In view of this fact, LZFC redirects our interest from absolute infinite cardinals to transitive models of ZFC which are now construed as analogues of inaccessible cardinals. Additional or stronger properties can  make models look like  stronger  large cardinals. For instance in \cite{Tz10} we  defined and studied   Mahlo and $\Pi_1^1$-indescribable models of ZFC, as analogues of Mahlo and weakly compact cardinals, respectively. By analogy we refer to such models as {\em large  models} of ZFC, and the corresponding properties as  {\em large model properties.} In this paper we continue the search for new kinds of large transitive  models of ZFC. So let us state  from the outset the following:

 \vskip 0.1in

{\bf Convention.} Throughout the paper ``model'' means ``transitive model of ZFC''.

\vskip 0.1in

Of course not all  large cardinal properties are expected to have sensible analogues for models.  For example such are the properties formulated in terms of ultrafilters.  In the absence of the Powerset axiom non-principal ultrafilters cannot be shown to exist even as   proper classes.  On the other hand, some properties formulated in terms of ultrafilters are equivalently formulated  over ZFC in terms of elementary embeddings, a  notion particularly fitted to a world with a plethora of local  models. So the properties formulated in the next section concern mostly existence of embeddings. The difference is that the embeddings employed in large cardinal definitions  are mappings $j:V\rightarrow W$ (where $W$ is an inner model), which are {\em internal}, i.e., definable in the universe $V$, while here  we deal with embeddings $j:M\rightarrow N$, between  models $M,N$ of ZFC,  which are in general {\em external} to both $M$ and $N$.

The main results of the paper are contained in section 3 (Theorem \ref{T:surprise}), section 5 (Theorem \ref{P:holds} and Proposition  \ref{P:allcountable}) and section 6 (Theorem \ref{T:main}).

The paper is organized as follows. In  section 2 we briefly review $\Pi_1^1$-indescribable models, which have been introduced in \cite{Tz10}, and give a simpler characterization of them in Lemma \ref{L:simplify}.

In section 3 we  introduce new stronger large model properties.  To every such property  $\phi(x)$ there naturally corresponds the localization principle $Loc^{\phi}({\rm ZFC})$ saying that ``every set belongs to a  model satisfying $\phi$''.  The theories
$${\rm LZFC}^\phi={\rm LZFC}+Loc^{\phi}({\rm ZFC})$$  extend  LZFC in  the same spirit as the theories $${\rm ZFC}+ \mbox{``there is a proper class of cardinals $\psi$''},$$
where $\psi$ is a large cardinal property,  extend  ZFC.

Firstly we consider  ``(elementarily) extendible'' models and show that every such model is $\Pi_1^1$-indescribable. Further, ``strongly extendible'' models are defined and  shown to be $\Sigma_2$-elementary submodels of the universe $V$. In particular if $sext(x)$ is the property of strong extendibility, the theory  ${\rm LZFC}^{sext}$ proves Powerset and $\Sigma_1$-Collection.
This  theory, as well as variants of it, are also largely  used  in section 6.

In section 4 we consider   elementary embeddings $j:M\rightarrow N$ between  models and their ``critical models''.   Stronger notions like ``strongly critical'', as well as  ``self-critical'' and ``strongly self-critical models'', are also introduced and some consequences of ZFC+$Loc^\phi({\rm ZFC})$, for various properties $\phi$, are proved.  Apparently critical models are the analogues of  measurable cardinals. However, due to the fact that the embeddings $j:M\rightarrow N$ are external, the consequences of their existence are much weaker than those of measurable cardinals. For example the  existence of strongly critical models alone does not seem to yield  $V\neq L$ over LZFC.

For that purpose in section 5 we introduce the Tall Model Axiom ($TMA$) saying, roughly, that for every ordinal $\kappa$ there is some ordinal $\alpha>\kappa$ such that there are models of arbitrarily big height which do not collapse $\alpha$ to $\kappa$. One of the main results of this section  (Theorem \ref{P:holds}) says that the theory LZFC+$TMA$+``there is a strongly critical model'' proves $V\neq L$. Further, the axiom of Greatest Cardinality ($GC$), saying that there is a (set of) greatest cardinality, is introduced. In the presence of $V=L$, this is equivalent to $\neg(TMA)$ over LZFC. Therefore LZFC+$TMA$+$GC$  yields $V\neq L$. However it is open whether LZFC+$TMA+GC$ is consistent.

In section 6 we examine the implications  of Vop\v{e}nka's Principle ($VP$) when added to LZFC. The key fact underlying these implications  is an old  ZFC result of  P. Vop\v{e}nka, A. Pultr and Z. Hedrl\'{i}n (abbreviated V-P-H) saying that for every set $A$ there is a relation $R\subset A\times A$ such that $(A,R)$ has no non-trivial endomorphism. The main result of section 6 is Theorem \ref{T:main} saying that  if $T$ is a theory such that ${\rm LZFC}\subseteq T$ and $T\vdash$ V-P-H, then ${\rm ZFC}\subseteq T+VP$, i.e., $T$ fully restores ZFC. An example of such a theory  is ${\rm LZFC}^{sext}$, as well as some weaker variants of it. Moreover the theories LZFC$^{sext}$+$VP$ and ZFC+$VP$ are identical.  It is not known whether  LZFC+$VP$ alone restores ZFC. However if LZFC+$VP$ proves either V-P-H or $Loc^{sext}({\rm ZFC})$, this is indeed the case. In connection to this it is shown that LZFC+$VP\vdash Loc^{ext}({\rm ZFC})$.

\section{$\Pi^1_1$ indescribable models}
Let us recall from \cite{Tz10} the definition of $\Pi_1^1$-indescribability.  Let ${\cal L}=\{\in\}$, let ${\cal L}_2$ be  ${\cal L}$ augmented with second order variables, and let ${\cal L}_2\cup\{{\bf S},{\bf c}_i\}$ be ${\cal L}_2$ augmented with   a unary predicate ${\bf S}(\cdot)$   intended to be  interpreted as a set $U$, and sufficient amount of  first-order constants ${\bf c}_i$. When $M$ is a model of ZFC, then the constants ${\bf c}_i$ are chosen to be names of  elements $c_i$ of $M$.  $Def(M)$ denotes the set of subsets of $M$ definable by formulas of ${\cal L}$. $Def(M)$ can be proved to exist in LZFC because it is absolute and can be constructed inside any model $N$ of ZFC that contains $M$ as a member.

\begin{Def} \label{D:transfer}
{\rm (LZFC)} {\em A  model $M\models{\rm ZFC}$ is said
to be} $\Pi^1_1$-indescribable {\em if for every $U\in Def(M)$ and every $\Pi^1_1$ sentence $\phi$ of ${\cal L}_2({\bf S},{\bf c}_i )$, if $(M,\in, U, Def(M))\models
\phi$, then there is a  model $N\in M$ such that $U\cap
N\in Def(N)$ and $(N,\in, \linebreak U\cap N, Def(N))\models
\phi$.}
\end{Def}

Actually in the above definition the set $U\cap N$ can be taken to be defined in $N$ by the same formula as $U$ in $M$, which  simplifies things considerably. Indeed, given a first order formula $\theta(x,\bar{y})$ without parameters, a model $M$ of ZFC, and $\bar{c}\in M$, let $\theta[M, \bf \bar{c}]$, or just $\theta[M]$ denote the set $\{x\in M:M\models\theta(x, \bf \bar{c})\}$. Then we have the following characterization of $\Pi^1_1$-indescribability (not contained in \cite{Tz10}).

\begin{Lem} \label{L:simplify}
A  model $M\models{\rm ZFC}$ is  $\Pi^1_1$-indescribable iff for every first-order formula $\theta(x,\bar{y})$ without parameters and every $\Pi^1_1$ sentence $\phi$, if ~~~~~~~~~~~~~~~~~~~~~~~~~~~~~~~~ $(M,\in, \theta[M], Def(M))\models \phi$, and ${\bar c}\in M$, then there is a model $N\in M$ such that $\bar{c}\in N$, $(N,\in, \theta[N], Def(N))\models\phi$,  and $\theta(x,\bf \bar{c})$ is absolute between $M$ and $N$ (i.e., $M\models\theta(x,{\bf \bar{c}})\Leftrightarrow  N\models \theta(x,\bf \bar{c})$, for every $x\in N$).
\end{Lem}

{\em Proof.} Suppose  $M$ is $\Pi^1_1$-indescribable. Let $\theta(x,\bar{y})$ be a first order formula, $\bar{c}\in M$,  $U=\theta[M]$, and  let  $\phi$ be a $\Pi^1_1$-formula of ${\cal L}_2({\bf S},{\bf c}_i )$ such that $(M,\in,U, Def(M))\models \phi$.  Let $\sigma:=\forall x({\bf S}(x)\leftrightarrow \theta(x,\bf \bar{c}))$. Then $(M,\in,U,Def(M))\models \phi \wedge \sigma$, and $\phi\wedge \sigma$ is $\Pi^1_1$.
By $\Pi^1_1$-indescribability there is a   model $N\in M$ such that $U\cap
N\in Def(N)$ and $(N,\in,U\cap N, Def(N))\models
\phi \wedge \sigma$. But $(N,\in,U\cap N, Def(N))\models
\sigma$ means that $U\cap N=\{x:N\models \theta(x,{\bf \bar{c}})\}=\theta[N]$,  thus $(N,\in,\theta[N], Def(N))\models
\phi$. Moreover, since  $\theta[N]=\theta[M]\cap N$, it follows that for every $x\in N$, $M\models\theta(x,{\bf \bar{c}})\Leftrightarrow  N\models \theta(x,{\bf \bar{c}})$, i.e., $\theta(x,\bf \bar{c})$ is absolute between  $M$ and $N$.

Conversely, suppose $M$ is a  model for which the assumption of the lemma holds. We show that $M$ is $\Pi^1_1$-indescribable. Let $U\in Def(M)$ and let $\phi$ be a $\Pi^1_1$ sentence such that $(M,\in,U,Def(M))\models \phi$. Let $U=\theta[M]$ for some $\theta$. By our assumption there is a model $N\in M$ such that $(N,\in,\theta[N], Def(N))\models \phi$, and  $\theta$ is absolute between $M$ and $N$. By the last  condition it follows that $\theta[N]=\theta[M]\cap N$. Hence $(N,\in,U\cap N, Def(N))\models\phi$. So $M$ is $\Pi^1_1$-indescribable. \telos

\vskip 0.2in

It was shown in  \cite[Proposition 5.5]{Tz10} that   $\Pi_1^1$-indescribability implies $\alpha$-Mahloness for every  model $M$ of ZFC. $\alpha$-Mahloness is an absolute (i.e., $\Delta_1$) property  in our formal language ${\cal L}=\{\in\}$ (see  \cite{Tz10}). The same is true for  the property of  $\Pi_1^1$-indescribability.
It can  be formalized by the  following absolute  formula $\pi_1^1ind(x)$:
$$\pi_1^1ind(x):=[x\models{\rm ZFC} \ \wedge \ (\forall y\in Def(x))(\forall \phi\in \Pi_1^1)(\exists z\in x)$$ $$(z\models{\rm ZFC} \wedge (x,\in, y,Def(x))\models \phi \Rightarrow y\cap z\in Def(z) \wedge (z,\in, y\cap z,Def(z))\models \phi)],$$
where $\Pi_1^1$ is the set of (codes of) $\Pi_1^1$-formulas of ${\cal L}_2({\bf S},{\bf c}_i )$.
Then for any  models $M,N$ such that $M\in N$, ``$M$ is  $\Pi_1^1$-indescribable'' (in $V$) iff $N\models \pi_1^1ind(M)$.

\section{Extendible and strongly extendible models}
In this section we introduce  properties   for    models of ZFC  stronger than $\Pi_1^1$-indescribability.

\begin{Def} \label{D:lee}
{\em  A  model $M$ is said to be}  elementarily  extendible, {\em or just} extendible, {\em if there is a  model $N$ such that $M\in N$ and $M\prec N$. }
\end{Def}

The next Lemma says that in the above definition the condition ``$M\in N$'' is redundant. This is a consequence of the fact  that for  models of ZFC, $M\prec N$ implies that $N$ is also an {\em end-extension} of $M$, i.e., the ordinals of $N$ extend those of $M$. Details are left to the reader.

\begin{Lem} \label{L:proper}
$M$ is  extendible  iff there is a model $N$ such that  $M\prec N$. In fact, $M\prec N\Rightarrow M\in N$.
\end{Lem}

The following induction scheme for the class $On$  is not part of LZFC.\footnote{See \cite[Lemma 2.11]{Tz10}, where  ${\rm Found}_{On}$ is shown to be equivalent to the $\in$-induction scheme
$$({\rm Found}^*) \hspace{.5\columnwidth minus .5\columnwidth}   \exists x\phi(x) \rightarrow \exists
x[\phi(x) \wedge \forall y\in x\neg\phi(y)]\hspace{.5\columnwidth minus .5\columnwidth} \llap{}$$
over LZFC.}
$$({\rm Found}_{On}) \hspace{.5\columnwidth minus .5\columnwidth}  \exists \alpha\in On \ \phi(\alpha) \rightarrow \exists
\alpha\in On[\phi(\alpha) \wedge \forall \beta<\alpha \neg\phi(\beta)]\hspace{.5\columnwidth minus .5\columnwidth} \llap{}$$

\begin{Prop} \label{P:why}
(i) In {\rm LZFC}: If $M$ is extendible, then $M$ is   $\Pi_1^1$-indescribable.

(ii) In ${\rm LZFC}+{\rm Found}_{On}$: The converse of (i) is false. I.e., if there is a $\Pi_1^1$-indescribable model, then there is one which is not  extendible.
\end{Prop}

{\em Proof.} (i) Suppose  $M$ is not $\Pi_1^1$-indescribable.
Then, by Lemma \ref{L:simplify},   there are a   $\Pi_1^1$ formula $\phi$ and a first-order formula  $\theta$ such that  $(M,\in, \theta[M],Def(M))\models\phi$, and for every  model $x\models {\rm ZFC}$ such that $x\in M$, either $\theta$ is not absolute between $M$ and $x$, or \linebreak $(x,\in, \theta[x], Def(x))\not\models\phi$. This means that  $(M,\in, \theta[M],Def(M))\models\phi$ and $M\models\psi_{\phi,\theta}$, where  $\psi_{\phi,\theta}$ is the first order formula
$$\psi_{\phi,\theta}:= \forall x[x\models{\rm ZFC}\wedge (\forall  y\in x)(\theta^x(y)\leftrightarrow \theta(y)) \rightarrow $$$$\hspace{1in}(x,\in, \theta[x],Def(x))\not\models\phi].$$
Suppose  that $M$ is  extendible, i.e.,   there is $N$ such that $M\in N$ and  $M\prec N$. Then  $N\models \psi_{\phi,\theta}$, i.e.,
$$N\models \forall x[x\models{\rm ZFC}\wedge (\forall  y\in x)(\theta^x(y)\leftrightarrow \theta(y)) \rightarrow $$
\begin{equation} \label{E:contra}
\hspace{1in}(x,\in, \theta[x],Def(x))\not\models\phi].
\end{equation}
But $M\in N$,  $M\models{\rm ZFC}$, $(\forall y\in M)(\theta^M(y)\leftrightarrow \theta^N(y))$, since $M\prec N$,   and   $(M,\in,\theta[M],Def(M))\models\phi$. This contradicts (\ref{E:contra}) and proves the claim.

(ii) Suppose (in ${\rm LZFC}+{\rm Found}_{On}$) that there are $\Pi_1^1$-indescribable models. By ${\rm Found}_{On}$  there is one such model $M$ of least height. Obviously $M$ cannot contain a $\Pi_1^1$-indescribable model. Let $\pi^1_1ind(x)$ be the formalization of the  $\Pi_1^1$-indescribability property.  $\pi^1_1ind(x)$ is absolute as noticed above.  Then $M$ is not  extendible. Assume the contrary. Then  there is a model $N$ such that $M\in N$ and $M\prec N$. Since $\pi^1_1ind(M)$ and $\pi^1_1ind(x)$ is absolute, it follows that $N\models \pi^1_1ind(M)$. Therefore $N\models \exists x(x\models{\rm ZFC} \wedge \pi^1_1ind(x))$. But then $M\models \exists x(x\models{\rm ZFC} \wedge \pi^1_1ind(x))$, which contradicts the fact that $M$ does not contain $\Pi^1_1$-indescribable models. \telos

\vskip 0.2in

As  an immediate  generalization of \ref{P:why} (ii) we have the following:

\begin{Prop} \label{P:minimal}
{\rm (${\rm LZFC}+{\rm Found}_{On}$)} If $\phi(x)$ is an absolute property about models and  $M$ is a  model of least height such that $\phi(M)$, then $M$ is not extendible.
\end{Prop}

In  contrast to  Mahloness and $\Pi^1_1$-indescribability, which are absolute properties, the formalization of  extendibility is given by the following  predicate:
$$ext(x):=[(x,\in)\models{\rm ZFC} \wedge (\exists y)((x,\in)\prec (y,\in))].$$

By analogy to large cardinal properties we may refer to  properties like $ext(x)$, $mahlo(x)$, $\pi^1_1in(x)$ (as well as to those that will be introduced later),  as {\em large model properties.} Analogously we may refer to   existence axioms  $\exists x \ ext(x)$, $\exists x \ mahlo(x)$,  etc,  as  {\em large model axioms.} Also for every large model property $\phi$ let
$$Loc^\phi({\rm ZFC}):=\forall x\exists y (x\in y \wedge y \ \mbox{transitive}  \wedge  \phi(y) \wedge   y\models {\rm ZFC}).$$
Namely, $Loc^\phi({\rm ZFC})$ says that every set $x$ belongs to a model satisfying the large model property $\phi$. So
we shall refer to $Loc^\phi({\rm ZFC})$ as {\em strong localization axioms,} since they strengthen  $Loc({\rm ZFC})$.
Moreover, for any such $\phi$ as above,
we shall denote by ${\rm LZFC}^\phi$ the theory resulting from LZFC if we replace the axiom $Loc({\rm ZFC})$ by $Loc^\phi({\rm ZFC})$.
That is:
$${\rm LZFC}^\phi:={\rm LZFC}+Loc^{\phi}({\rm ZFC})={\rm BST}+Loc^{\phi}({\rm ZFC}).$$
The theories ${\rm LZFC}^\phi$ are extensions of LZFC, pretty analogous to the extensions of ZFC of the form
$${\rm ZFC}+ \mbox{``there is a proper class of cardinals $\psi$''},$$
for some  large cardinal property $\psi$. Later in this section as well as in section 6, we  shall largely  work in the theory ${\rm LZFC}^{sext}$, where  $sext(x)$ is the property of strong extendibility (see Definition \ref{D:sext} below), and some variants of it.

The next  result of ZFC  is a slight generalization  of  \cite[Theorem 8.1]{MV59}) and   will be needed in sections 5 and 6. The proof is similar to the standard one and left to the reader.

\begin{Lem} \label{L:montague}
{\rm (ZFC)} If there are ordinals $\alpha<\beta$ such that  $V_\alpha\prec V_\beta$, then  $V_\alpha, V_\beta$ are models of {\rm ZFC}. More generally, if $\alpha>\omega$, $(V_\alpha,\in)\prec (A,\in)$ and $V_\alpha\in A$,  for some transitive set $A$, then $V_\alpha\models {\rm ZFC}$.
\end{Lem}

For every axiomatized set theory $T$ (where usually $T\supseteq {\rm ZFC}$), let $Loc(T)$ denote the corresponding localization principle for $T$, i.e.,
$$Loc(T):=\forall x\exists
y(x\in y \ \wedge \ Tr(y) \ \wedge \ (y,\in)\models T).$$

We have seen in \cite{Tz10} that every Mahlo model of ZFC satisfies $Loc_n({\rm ZFC})$, for every $n\in\omega$, where these  properties are defined inductively by the clauses: $Loc_{0}({\rm ZFC})=Loc({\rm ZFC})$ and $Loc_{n+1}({\rm ZFC})=Loc({\rm ZFC}+Loc_n({\rm ZFC}))$. More generally, in \cite{Tz10} we considered strengthenings of $Loc({\rm ZFC})$, of the form $Loc({\rm ZFC}+\phi)$, where $\phi$ is a new  axiom added to ZFC.

The sentences $Loc_n({\rm ZFC})$ are  estimates of how much stronger  a Mahlo model is compared to  an ordinary model of ZFC. The question is whether  there are  analogous estimates for $\Pi_1^1$-indescribable and elementarily extendible  models. We see below that the answer is yes.

We have seen that the properties of Mahloness, $\Pi_1^1$-indescribability and elementary extendibility are in increasing strength (for the last two, when working in LZFC+$Found_{On}$). But a natural question is: In what sense are the models of one of these classes {\em internally stronger} than  those of another. For instance, what properties does a $\Pi_1^1$-indescribable model satisfy which   a Mahlo model does  not?  The next two lemmas reasonably justify this internal increasing strength.
Recall that $mahlo_\alpha(x)$ abbreviates  the property $mahlo(\alpha,x)$ defined in \cite{Tz10}.

\begin{Lem} \label{L:isittrue}
If $M$ is $\Pi_1^1$-indescribable, then $M\models Loc^{mahlo_\alpha}({\rm ZFC})$, for all $\alpha\in On\cap M$.
\end{Lem}

{\em Proof.} Let $M$ be $\Pi_1^1$-indescribable. Fix some $a\in M$ and let $\alpha\in On^M$. We have to show that the set $X=\{x\in M:a\in x \wedge x\models{\rm ZFC} \wedge mahlo(\alpha,x)\}\neq\emptyset$. Let $Y=\{x\in M:a\in x \}$. Obviously, $Y$ is a club of $M$. By Proposition 5.5 of \cite{Tz10}, since $M$ is $\Pi_1^1$-indescribable, $M$ is $(\alpha+1)$-Mahlo. This means that  the set $Z=\{x\in M:x\models {\rm ZFC} \wedge mahlo(\alpha,x)\}$ is stationary in $M$. Hence $Z\cap Y\neq \emptyset$. But $Z\cap Y=X$, therefore $X\neq \emptyset$. \telos

\begin{Lem} \label{L:isittrue1}
If $M$ is  extendible then   $M\models Loc^{\pi^1_1ind}({\rm ZFC})$.
\end{Lem}

{\em Proof.} Let $M$ be  extendible and let $a\in M$. We have to show that $X=\{x\in M:a\in x\wedge x\models{\rm ZFC} \wedge \pi^1_1ind(x)\}\neq \emptyset$.  Suppose that $X=\emptyset$. Since $\pi^1_1ind(x)$ and ``$x\models{\rm ZFC}$'' are  absolute, this is equivalent to
\begin{equation} \label{E:exo}
M\models \forall x(a\in x \wedge x\models{\rm ZFC} \rightarrow  \neg \pi^1_1ind(x)).
\end{equation}
By our assumption there is a model $N$ such $M\in M$ and $M\prec N$. Then, by (\ref{E:exo}),
\begin{equation} \label{E:exo1}
N\models \forall x(a\in x \wedge x\models{\rm ZFC} \rightarrow  \neg \pi^1_1ind(x)).
\end{equation}
But  by Proposition \ref{P:why}, $M$ is $\Pi^1_1$-indescribable, i.e., $\pi^1_1ind(M)$, the latter property being absolute. So $a\in M$, $M\in N$, $M\models{\rm ZFC}$ and $N\models\pi^1_1ind(M)$, which contradicts (\ref{E:exo1}). \telos

\vskip 0.2in

The following  strengthening of  extendibility  turns out to be an interesting and powerful  property.

\begin{Def} \label{D:sext}
{\em A model $M$ is} strongly  extendible {\em if for every $x$ there is a model $N$ such that  $x\in N$ and $M\prec N$.}
\end{Def}
Formally the property of strong extendibility is written:
$$sext(x):=[(x,\in)\models{\rm ZFC} \wedge (\forall y)(\exists z)(y\in z \wedge (x,\in)\prec (z,\in))].$$

Given a model $M$ let us write $M\prec_{\Sigma_n} V$ if every $\Sigma_n$ sentence  $\phi$ with parameters from $M$ is absolute for $M$, i.e., $M\models \phi$ iff $V\models \phi$. [A $\Sigma_n$ formula is one of the form $(\exists \overline{x}_1)(\forall \overline{x}_2)\cdots (Q \overline{x}_n)\psi$, with   $\psi$ bounded, where $\overline{x}_i$ are tuples of variables. Similarly for a $\Pi_n$ formula.]   Obviously $M\prec_{\Sigma_n} V$ iff $M\prec_{\Pi_n} V$.

The important feature of strongly extendible  models is the following.

\begin{Lem} \label{L:abssigma}
Let  $M$ be strongly extendible.  Then:

(i) $M\prec_{\Sigma_1} V$. Therefore, in ${\rm LZFC}^{sext}$  every set  belongs to some model $M\prec_{\Sigma_1} V$.

(ii) Actually $M\prec_{\Sigma_2}V$.
\end{Lem}

{\em Proof.} Fix a strongly extendible  model $M$.

(i) Let $\phi=(\exists \overline{x})\psi(\overline{x},\overline{c})$ be a $\Sigma_1$ sentence, with $\psi$ bounded and  $\overline{c}\in M$. Then clearly $M\models (\exists \overline{x})\psi(\overline{x},\overline{c})$ implies $V\models(\exists \overline{x})\psi(\overline{x},\overline{c})$. Conversely, assume $V\models(\exists \overline{x})\psi(\overline{x},\overline{c})$, and let $V\models\psi(\overline{a},\overline{c})$ for some tuple $\overline{a}$. By strong extendibility  there is a model $N$ such that $M\prec N$ and $\overline{a}\in N$. By the absoluteness of $\psi$, $N\models \psi(\overline{a},\overline{c})$, which implies $N\models (\exists \overline{x}) \psi(\overline{x},\overline{c})$. Since $M\prec N$, it follows that $M\models (\exists \overline{x}) \psi(\overline{x},\overline{c})$.\footnote{Actually the proof of this clause does  not require the full strength of the condition of strong extendibility. It is easy to see that the following weaker condition for $M$, called ``$\Sigma_1$-strong extendibility'',  suffices: For every $x$ there is a model $N$ such that $x\in N$ and $M\prec_{\Sigma_1}N$. See section 6 for more on this property. \label{f:4}}

(ii) Let $\phi=(\exists \overline{x})(\forall \overline{y})\psi(\overline{x},\overline{y},\overline{c})$ be a $\Sigma_2$ sentence with $\psi$ bounded and  $\overline{c}\in M$, and let $M\models (\exists \overline{x})(\forall \overline{y})\psi(\overline{x},\overline{y},\overline{c})$. Then for some $\overline{a}\in M$, $M\models (\forall \overline{y})\psi(\overline{a},\overline{y},\overline{c})$. The sentence $(\forall \overline{y})\psi(\overline{a},\overline{y},\overline{c})$ is $\Pi_1$, so by (i) above, $V\models (\forall \overline{y})\psi(\overline{a},\overline{y},\overline{c})$. Therefore $V\models (\exists \overline{x})(\forall \overline{y})\psi(\overline{x},\overline{y},\overline{c})$. Conversely, assume $V\models (\exists \overline{x})(\forall \overline{y})\psi(\overline{x},\overline{y},\overline{c})$. Then $V\models (\forall \overline{y})\psi(\overline{a},\overline{y},\overline{c})$ for some $\overline{a}$. By the strong extendibility of $M$ there is a model $N$ such that $\overline{a}\in N$ and $M\prec N$. Since $N$ is a transitive submodel of $V$,   $(\forall \overline{y})\psi(\overline{a},\overline{y},\overline{c})$ is $\Pi_1$, and the latter is true  in $V$, it follows that  $N\models (\forall \overline{y})\psi(\overline{a},\overline{y},\overline{c})$. Consequently, $N\models (\exists \overline{x})(\forall \overline{y})\psi(\overline{x},\overline{y},\overline{c})$. It follows that  $M\models (\exists \overline{x})(\forall \overline{y})\psi(\overline{x},\overline{y},\overline{c})$. \telos

\vskip 0.2in

Already $Loc^{sext}({\rm ZFC})$ partly restores ZFC. Namely:

\begin{Thm} \label{T:surprise}
(i) ${\rm LZFC}^{sext}\vdash {\rm Powerset}$.

(ii) ${\rm LZFC}^{sext}\vdash \Sigma_1$-{\rm Collection}.
\end{Thm}

{\em Proof.} (i) Let $a$ be a set and let $M$ be a strongly extendible model such that $a\in M$. Let $b={\cal P}^M(a)$ be the powerset of $a$ in $M$, i.e., $M\models b={\cal P}(a)$.  The predicate $y={\cal P}(x)$ is $\Pi_1$, so by \ref{L:abssigma} (i),  $V\models b={\cal P}(a)$, so $b$ is the absolute powerset of $a$.

(ii) Let $\phi(x,y,\overline{c})=(\exists \overline{z})\psi(x,y,\overline{z},\overline{c})$ be a $\Sigma_1$ formula with parameters $\overline{c}$, let $a$ be a set and let
$(\forall x\in a)(\exists y)\phi(x,y,\overline{c})$ be true  (in $V$). We have to show that there is a set $b$ such that $(\forall x\in a)(\exists y\in b)\phi(x,y,\overline{c})$. In ${\rm LZFC}^{sext}$ we can pick a model $M$ such that $a,\overline{c}\in M$ and $M$ is strongly extendible. By \ref{L:abssigma} (ii),  $M\prec_{\Pi_2} V$ and by assumption $V\models (\forall x\in a)(\exists y)\phi(x,y,\overline{c})$, or $$V\models (\forall x\in a)(\exists y)(\exists \overline{z})\psi(x,y,\overline{z},\overline{c}).$$  The last formula is $\Pi_2$, so $M\models (\forall x\in a)(\exists y)(\exists \overline{z})\psi(x,y,\overline{z},\overline{c})$. Since $M$ satisfies Collection, there is a $b\in M$ such that $$M\models (\forall x\in a)(\exists y\in b)(\exists \overline{z})\psi(x,y,\overline{z},\overline{c}).$$ By $M\prec_{\Pi_2} V$ again, $$V\models (\forall x\in a)(\exists y\in b)(\exists \overline{z})\psi(x,y,\overline{z},\overline{c}),$$ i.e., $V\models (\forall x\in a)(\exists y\in b)\phi(x,y,\overline{c})$, as required. \telos

\vskip 0.2in

Further, sets  $V_\alpha={\cal P}^\alpha(\emptyset)$ can be defined  in ${\rm LZFC}^{sext}$ for every $\alpha\in On$ without the help of induction along $On$, i.e., ${\rm Found}_{On}$.  Moreover, as in the case of ZFC, $V_\alpha$ are set approximations of $V$  and include all strongly extendible models.

\begin{Prop} \label{P:rankinside}
In  ${\rm LZFC}^{sext}$ the following hold:

(i) Given any  specific ordinal $\alpha$,  for all $\beta\leq \alpha$, there are transitive sets $V_\beta={\cal P}^\beta(\emptyset)$ with the usual {\rm ZFC} properties, namely,   $V_{\beta+1}={\cal P}(V_\beta)$ and $V_\beta=\bigcup_{\gamma<\beta}V_\gamma$, for limit $\beta$. Moreover, for every strongly extendible model $M$ such that  $\alpha\in M$, $M_\beta:=V_\beta^M=V_\beta$, for all $\beta\leq \alpha$.

(ii) The predicate ``$x=V_\alpha$'' is definable, therefore so is the class $\{V_\alpha:\alpha\in On\}$.

(iii)  For every strongly extendible $M$, there is $\alpha$ such that $M=V_\alpha$.

(iv) $V=\bigcup_\alpha V_\alpha$. Moreover there are arbitrarily large $\alpha$ such that $V_\alpha\prec_{\Sigma_2}V$.
\end{Prop}

{\em Proof.} (i) Given an ordinal $\alpha$ pick, by $Loc^{sext}({\rm ZFC})$, a strongly extendible $M$ such that $\alpha\in M$. By Theorem \ref{T:surprise} (i) and induction inside $M$,
$$M_\beta=V_\beta^M=({\cal P}^\beta(\emptyset))^M={\cal P}^\beta(\emptyset)=V_\beta$$ for all $\beta\leq \alpha$.

(ii) $V_\alpha$ can be defined by induction up to $\alpha$  inside any strongly extendible model containing $\alpha$, the steps of which will be absolute because of \ref{T:surprise} (i). Namely ``$x=V_\alpha$'' is written:
$$(\exists f)[dom(f)=\alpha+1 \ \wedge \ (\forall \beta< \alpha)(f(\beta+1)={\cal P}(f(\beta))\ \wedge$$ $$(\forall \beta\leq \alpha)(\beta \ \mbox{limit} \ \Rightarrow f(\beta)=\bigcup_{\gamma<\beta}f(\gamma)) \ \wedge \ f(\alpha)=x].$$

(iii) Let $M$ be strongly extendible. We can show that $\alpha=\min\{\beta:\beta\notin M\}$, the height of $M$, exists  without invoking ${\rm Found}_{On}$. Indeed, given $M$ just take a model $N$, by $Loc({\rm ZFC})$, such that $M\in N$. Then clearly $\alpha=\min\{\beta:\beta\notin M\}$ exists in $N$. By (i) above  $V_\alpha$ exists too and   for every $\beta\in M$, $V_\beta=M_\beta$, therefore $$V_\alpha=\bigcup_{\beta<\alpha}V_\beta=
\bigcup_{\beta<\alpha}M_\beta=M.$$

(iv) Given $x$, by $Loc^{sext}({\rm ZFC})$ there is a strongly extendible $M$ such that $x\in M$. But by (iii) above, $M=V_\alpha$ for some $\alpha$, thus $\bigcup_\alpha V_\alpha=V$. Also, by \ref{L:abssigma} (ii), $V_\alpha\prec_{\Sigma_2} V$ for every strongly extendible $V_\alpha$.  \telos

\vskip 0.2in

What is the  consistency strength of $Loc^{sext}({\rm ZFC})$?  Surprisingly enough, this is no higher than the existence of a strongly inaccessible cardinal. Namely the following holds.

\begin{Prop} \label{P:nostrong}
\rm{(ZFC)} If $\kappa$ is strongly inaccessible, then $V_\kappa\models Loc^{sext}({\rm ZFC})$.
\end{Prop}

{\em Proof.} It is well-known  that if $\kappa$ is strongly inaccessible, then  the set $C=\{\alpha<\kappa:(V_\alpha,\in)\prec (V_\kappa,\in)\}$ is closed unbounded (see \cite[p. 171]{Je03}). Given $\alpha<\beta$ in $C$, we have  $V_\alpha\prec V_\kappa$ and $V_\beta\prec V_\kappa$, and therefore $V_\alpha\prec V_\beta$. Since $C$ is cofinal in $\kappa$, it follows immediately that for every $\alpha\in C$, $V_\alpha$ is strongly extendible in $V_\kappa$.  Also every $x\in V_\kappa$ belongs to  $V_\alpha$, for some $\alpha\in C$, so $V_\kappa\models Loc^{sext}({\rm ZFC})$. \telos

\vskip 0.2in

One further question: For which cardinals $\kappa$ of ZFC is $V_\kappa$ strongly extendible?  Note that if in the proof of \ref{P:nostrong} $\kappa$ is the {\em least} strongly inaccessible cardinal of the universe, then for each  $\alpha\in C$, $V_\alpha$ is a strongly extendible model, and moreover $V_\alpha\prec V_\kappa$, without $\alpha$ being strongly inaccessible. Therefore strong extendibility of $V_\kappa$, in the context of ZFC, does not presume even the inaccessibility of $\kappa$. However if we need a {\em sufficient} condition about $\kappa$ in order for $V_\kappa$ to be strongly extendible, we must go up to strong cardinals. Recall the following.

\begin{Def} \label{D:strongc}
{\rm (ZFC)} {\em A cardinal  $\kappa$ is} $\alpha$-strong {\em if there is an elementary embedding  $j:V\rightarrow W$ such that $\kappa={\rm crit}(j)$, $j(\kappa)>\alpha$ and $V_{\kappa+\alpha}\subseteq W$. $\kappa$ is} strong {\em if it is $\alpha$-strong for every $\alpha$. }
\end{Def}

\begin{Lem} \label{L:strong}
{\rm (ZFC)} If $\kappa$ is a strong cardinal then $V_\kappa$ is a strongly extendible model.
\end{Lem}

{\em Proof.} Let $\kappa$ be strong and $x$ be a set. We have to find a model  $M$ such that $V_\kappa\prec M$ and $x\in M$. Let $\alpha>\kappa$ be such that $x\in V_\alpha$ and $\kappa+\alpha=\alpha$. Since $\kappa$ is $\alpha$-strong,  there is $j:V\rightarrow W$ with ${\rm crit}(j)=\kappa$, $j(\kappa)>\alpha$  and $V_\alpha=V_{\kappa+\alpha}\subseteq W$. Therefore $V_\beta=W_\beta$ for every $\beta\leq \alpha$. Since ${\rm crit}(j)=\kappa$,  $V_\kappa\prec j(V_\kappa)=W_{j(\kappa)}$. Also $V_\alpha=W_\alpha\subset W_{j(\kappa)}$, thus $x\in W_{j(\kappa)}$. So for $M=W_{j(\kappa)}$, $M$ is as required. \telos

\section{Elementary embeddings and critical models}

\subsection{Elementary embeddings}

In this section we define models of ZFC with properties analogous to those of some large cardinals, employing  elementary embeddings.  Let $M$, $N$ be  models of ZFC and let  $j:M\rightarrow N$ be an elementary embedding. If $j$ is non-trivial, then   there is  a  critical ordinal ${\rm crit}(j)$ for $j$, i.e., a least ordinal $\alpha\in M$ such that $j(\alpha)>\alpha$. Moreover we can see that ${\rm crit}(j)$ is inaccessible in $M$. (This is essentially proved in \cite{Co06}.) More generally, let us call a set $x$ {\em critical} for $j$, if $j\restr x=id$ while $j(x)\neq x$. Let ${\rm Crit}(j)$ be the set of critical sets of $j$. In particular  ${\rm crit}(j)\in {\rm Crit}(j)$. In view of the axiom $Loc({\rm ZFC})$, given $M,N$ and  $j:M\rightarrow N$,  ${\rm Crit}(j)$ is a set in LZFC  because it can be (absolutely) defined in every model $K$ of ZFC  such that $\{M,N,j\}\subset K$.  The following two lemmas contain some basic facts about  elementary embeddings between arbitrary  models of ZFC. The proofs are mostly standard and scattered (usually in the form of exercises) in \cite{Je03}, \cite{Ka97}, etc. For the reader's convenience we outline  here some of them.

\begin{Lem} \label{L:inac}
{\rm (LZFC)} Let $j:M\rightarrow N$ be a nontrivial elementary embedding between  models of {\rm ZFC}. Then:

(i)  ${\rm crit}(j)$ exists.

(ii) If $\kappa={\rm crit}(j)$, then $j\restr M_\kappa=id$ (where as usual $M_\kappa=V_\kappa^M$). Therefore  $M_\kappa\in {\rm Crit}(j)$.  In particular, for every $A\in M_\kappa$, ${\cal P}(A)^M={\cal P}(A)^N$. Conversely, if $K\in {\rm Crit}(j)$ is a model,  then $ht(K)=\kappa$.

(iii) If $\kappa={\rm crit}(j)$, then $M_\kappa$ is the greatest transitive critical set of $M$, in fact if $x$ is a transitive set of $M$ and $j\restr x=id$, then $x\subseteq M_\kappa$.

(iv) If $\kappa={\rm crit}(j)$, then ${\cal P}(\kappa)^M\subseteq {\cal P}(\kappa)^N$.

(v) If $\kappa={\rm crit}(j)$, then $M\models$ ``$\kappa$ is strongly inaccessible''. Therefore $M_\kappa\models{\rm ZFC}$. If in addition $N=M$, then $M\models ``\mbox{$\kappa$ is $n$-ineffable}$'', for every $n\in\omega$.\footnote{A cardinal $\lambda$ is {\em $n$-ineffable} if for every partition $f:[\lambda]^{n+1}\rightarrow \{0,1\}$ there is a stationary homogeneous set $H\subseteq \lambda$.}\label{F:2}

(vi) For every first-order structure  $A\in M$, $j\restr A:A\rightarrow j(A)$ is an elementary embedding.  If in addition,  $A\in {\rm Crit}(j)$, then  $A\prec j(A)$. In particular, for the critical $\kappa$,  $M_\kappa\prec j(M_\kappa)$.
\end{Lem}

{\em Proof.} (i) For ZFC, the proof that $j$ has a critical point is given in \cite[Prop. 5.1 (b)]{Ka97} for the case $M\models {\rm AC}$. To show this in LZFC, given $M,N,j$, simply pick by  $Loc({\rm ZFC})$ a model  $K\models {\rm ZFC}$ such that $M,N,j\in K$ and work in $K$.

Clauses (ii)-(iv) are left to the reader.

(v) Let $\kappa={\rm crit}(j)$.
In \cite[p. 338]{Co06} it is proved (in ZFC) that if $V$ is a model of ZFC and $j:V\rightarrow V$ is an elementary embedding,  in general  external with respect to $V$, and ${\rm crit}(j)=\kappa$,  then $\kappa$ is inaccessible. The same proof actually works  for  general  elementary embeddings $j:M\rightarrow N$. Indeed, let $j:M\rightarrow N$ be a nontrivial elementary embedding and let  $\kappa={\rm crit}(j)$. We show first that $M\models ``\mbox{$\kappa$ is regular}$''. Obviously $\kappa$ is a limit ordinal. Assume on the contrary that there is $\alpha<\kappa$ and a cofinal $f:\alpha\rightarrow \kappa$, $f\in M$. Then $dom(f)=\alpha$, therefore $dom(j(f))=j(\alpha)=\alpha$. Moreover, for every $\beta<\alpha$, if $f(\beta)=\gamma$, then $\gamma<\kappa$, so $j(\beta)=\beta$ and $j(\gamma)=\gamma$. Therefore, $j(f)(\beta)=j(f)(j(\beta))=j(f(\beta)=j(\gamma)=\gamma=f(\beta)$. Thus $dom(j(f))=dom(f)$ and $j(f)(\beta)=f(\beta)$ for every $\beta\in dom(f)$, so $j(f)=f$. But since $f''\alpha$ is cofinal in $\kappa$, $j(f)''\alpha=f''\alpha$ should  be cofinal in $j(\kappa)>\kappa$, a contradiction.

Next suppose $\kappa$ is not a strong limit in $M$, i.e., there is $\alpha<\kappa$ such that $M\models |{\cal P}(\alpha)|\geq \kappa$. Now since $\alpha<\kappa$, by  (ii) above, ${\cal P}(\alpha)^M={\cal P}(\alpha)^N$.  Let $g\in M$ be a surjection $g:{\cal P}(\alpha)^M\rightarrow \kappa$.
Since $j({\cal P}(\alpha)^M)={\cal P}(\alpha)^N={\cal P}(\alpha)^M$,  it follows as before that $j(g)=g$. But then $j(g)=g$ must be a surjection of ${\cal P}(\alpha)^N$ onto $j(\kappa)>\kappa$, a contradiction. This completes the proof that $\kappa$ is strongly inaccessible in $M$. Therefore $M\models ``\mbox{$V_\kappa$ is a model of {\rm ZFC}}$'', which means that $M_\kappa\models{\rm ZFC}$. Finally, if $N=M$, then the proof of Theorem 2.12 of \cite{Co06} works also here, showing that the critical $\kappa$ is $n$-ineffable for every $n\in\omega$.

(vi)   For every formula $\phi(x_1,\ldots,x_n)$  of the language of $A$, and every $a_1,\ldots,a_n\in A$, in view of the absoluteness of satisfaction we have:
$$A\models\phi(a_1,\ldots,a_n)\Leftrightarrow
M\models(A\models\phi(a_1,\ldots,a_n))\Leftrightarrow$$
$$N\models(j(A)\models\phi(j(a_1),\ldots,j(a_n))\Leftrightarrow j(A)
\models\phi(j(a_1),\ldots,j(a_n)),$$
which means that  $j\restr A:A\rightarrow  j(A)$ is an elementary embedding. If  $A\in {\rm Crit}(j)$, then $j(a_i)=a_i$ for $a_i\in A$, while $j(A)\neq A$. Thus $A\models\phi(a_1,\ldots,a_n)\Leftrightarrow
j(A)\models\phi(a_1,\ldots,a_n)$, i.e., $A\prec j(A)$. \telos

\begin{Def} \label{D:embed}
{\em A  model $M$ of {\rm ZFC} is said to be} elementarily embeddable {\em or just} embeddable, {\em if there is a non-trivial elementary embedding $j:M\rightarrow N$ for some $N$. $M$ is said to be} inner-embeddable {\em if there is an elementary embedding $j:M\rightarrow N\subseteq M$. Finally, $M$ is} self-embeddable {\em if there is a (non-trivial) elementary embedding $j:M\rightarrow M$. }
\end{Def}

The relationship between elementary embeddability and elementary extendibility (considered in the previous section) is not quite clear. One might say that extendibility is stronger since if $M$ is extendible and $M\prec N$, then $M$ is also embeddable in $N$ with respect to the  trivial embedding $id$. But the essence of embeddability is exactly the existence of a {\em non-trivial} embedding which gives rise to critical points.  The two notions are probably incomparable.

As we have seen in \ref{L:isittrue1}, if $M$ is elementarily extendible  then   $M\models Loc^{\pi^1_1ind}({\rm ZFC})$. But if $M$ is embeddable or even self-embeddable, $M$ is unlikely to  satisfy even   $Loc({\rm ZFC})$. However the following   holds.

\begin{Lem} \label{L:infinac}
(i) If $M$ is embeddable, then $M\models$``there is a  model $x$ such that $x\models{\rm ZFC}+Loc({\rm ZFC})$''.

(ii) If $M$ is self-embeddable and $j:M\rightarrow M$ is an embedding with $\kappa={\rm crit}(j)$, then the cardinals $j^n(\kappa)$ are $n$-ineffable in $M$.
\end{Lem}

{\em Proof.} Let $j:M\rightarrow N$ be an elementary embedding. By \ref{L:inac} (v), if $\kappa={\rm crit}(j)$, $M\models$ ``$\kappa$ is strongly inaccessible''. But then  $M\models (V_\kappa\models {\rm  ZFC}+Loc({\rm ZFC}))$. That is, the model $x\in M$ in question  such that $x \models {\rm  ZFC}+Loc({\rm ZFC})$ is $M_\kappa$.

(ii) By \ref{L:inac} (v). \telos

\vskip 0.2in

Elementary embeddability is formalized by the following predicate:
$$emb(x):=[x\models{\rm ZFC} \ \wedge \ (\exists y)(\exists j)
(y\models{\rm ZFC} \ \wedge \  j:x\rightarrow  y \ \mbox{is el. emb.})].$$

\subsection{Critical and strongly critical models}

\begin{Def} \label{D:crit}
{\em $M$ is said to be} critical {\em if there are models $N,K$ such that $M\in N$ and an elementary embedding $j:N\rightarrow K$ such that $M\in {\rm Crit}(j)$. }
\end{Def}

The property is formalized also be the  $\Sigma_1^{\rm ZFC}$ predicate:
$$crit(x):=[x\models{\rm ZFC} \ \wedge \ (\exists y)(\exists z)(\exists j)$$
$$(x\in y \ \wedge \ y\models{\rm ZFC} \ \wedge \ z\models{\rm ZFC}\ \wedge \  j:y\rightarrow z \ \mbox{is el. emb} \wedge \ x\in {\rm Crit}(j))].$$

Critical models are  roughly analogues of  measurable cardinals of ZFC. And their consistency strength is no greater than that of the latter.

\begin{Lem} \label{L:measurec}
{\rm(ZFC)} If $\kappa$ is a measurable cardinal, then $V_\kappa$ is  critical.
\end{Lem}

{\em Proof.}  Let $\kappa$ be a measurable cardinal in {\rm ZFC}, and let $j:V\rightarrow W$ be an elementary embedding (where $W\subseteq V)$), with ${\rm crit}(j)=\kappa$.  Then clearly $j\restr V_\kappa=id$, while $j(V_\kappa)\neq V_\kappa$. So it suffices to find models $M, N$ such that $j\restr M:M\rightarrow N$ and $V_\kappa\in M$. Set   $M=j(V_\kappa)=W_{j(\kappa)}$,  and $N=j(M)$. Then, clearly $M,N$ are models of ZFC since $V_\kappa$ is so, $V_\kappa=W_\kappa\in W_{j(\kappa)}=M$, and if $j'=j\restr M$, then $j':M\rightarrow N$ is an elementary embedding with $V_\kappa\in {\rm Crit}(j')$.  \telos

\vskip 0.2in

\begin{Lem} \label{L:crit}
{\rm(LZFC)} If  $M$ is  critical, then   $M$ is extendible. Therefore over {\rm LZFC}, $Loc^{crit}({\rm ZFC})\Rightarrow Loc^{ext}({\rm ZFC})$.
\end{Lem}

{\em Proof.}  Let $M$ be  critical in LZFC. Then there are $N,K,j$ such that $M\in N$, $j:N\rightarrow K$ is an elementary embedding and $M\in {\rm Crit}(j)$. By Lemma \ref{L:inac} (vi), $M\prec j(M)$. Therefore $M$ is  extendible. \telos

\vskip 0.2in

The following strengthening of criticalness is reasonable:

\begin{Def} \label{D:stongrcrit}
{\em Let $M$ be a model and $x$ be a set.  $M$ is said to be} $x$-strongly critical {\em if there are models $N,K$ and  an elementary embedding $j:N\rightarrow K$, such that $\{M,x\}\subset N$ and  $M\in {\rm Crit}(j)$. $M$ is said to be} strongly critical {\em if it is $x$-strongly critical for every $x$.}
\end{Def}

Here is a simpler characterization of strong criticalness.

\begin{Lem} \label{L:chara}
{\rm (LZFC)} $M$ is strongly critical iff for every model $N$ such that $M\in N$, there is a model $K$ and an elementary embedding $j:N\rightarrow K$ such that $M\in {\rm Crit}(j)$.
\end{Lem}

{\em Proof.} The condition is necessary. Let $M$ be strongly critical and let $N$ be a model such $M\in N$. By definition $M$ is $N$-critical, i.e., there are models $R,S$ with $\{M,N\}\subset R$ and an elementary embedding $j:R\rightarrow S$ such that $M\in {\rm Crit}(j)$. Let $j'=j\restr N$ and $K=j(N)$. Then, according to Lemma \ref{L:inac} (vi),  $j':N\rightarrow K$ is an elementary embedding with $M\in {\rm Crit}(j')$. Conversely, suppose the condition holds for $M$ and let $x$ be a set. By $Loc({\rm ZFC})$ there is a model $N$ of ZFC such that
$\{ M,x\}\subset N$. By our condition there are $K$ and an elementary embedding $j:N\rightarrow K$ with $M\in {\rm Crit}(j)$. Thus $M$ is $x$-strongly critical for every $x$, and therefore is  strongly critical. \telos.

\begin{Prop} \label{P:Lalpha}
{\rm (LZFC)} Let  $M$ be a strongly  critical model of height $\kappa$. Then:

(i)  For every ordinal $\alpha>\kappa$, $L_\alpha$ is  embeddable into some $L_\beta$ with critical point $\kappa$. In particular  there is a cofinal class $C\subseteq On$, such that for every $\alpha\in C$,  $L_\alpha$ is a  model of {\rm ZFC} and $L_\alpha$ is embeddable into some $L_\beta$ with critical point $\kappa$.

(ii) If in addition $V=L$, then there is a  cofinal class $C\subseteq On$, such that for every $\alpha\in C$,  $L_\alpha$ is a  model of {\rm ZFC}, $M\in L_\alpha$ and $L_\alpha$ is embeddable into some $L_\beta$ with critical model $M$.
\end{Prop}

{\em Proof.} (i) Suppose  $M$ is a strongly  critical model of height $\kappa$. Let $\alpha>\kappa$. By strong criticalness, there are models $N,K$ such that $\{\alpha, M\}\subset N$ and an elementary embedding $j:N\rightarrow K$ such that ${\rm Crit}(j)=M$;  hence ${\rm crit}(j)=\kappa$. Then $L_\alpha\in N$. If $j'=j\restr L_\alpha$, then by Lemma  \ref{L:inac} (vi), $j':L_\alpha\rightarrow L_{j(\alpha)}$ is a nontrivial elementary embedding with ${\rm crit}(j)=\kappa$. As for the second claim, let $$C=\{\alpha\in On: \alpha>\kappa \wedge L_\alpha\models {\rm ZFC}\}.$$
By $Loc({\rm ZFC})$, $C$ is cofinal in $On$. By the preceding argument, for every $\alpha\in C$, there is a model $N$ such that $\alpha\in N$ and there is an embedding $j:N\rightarrow K$ such that $\kappa={\rm crit}(j)$. Then $j\restr L_\alpha:L_\alpha\rightarrow L_{j(\alpha)}$ has also critical point $\kappa$.

(ii) If $V=L$, then $M\in L$ and if we set
$$C=\{\alpha\in On: M\in L_\alpha  \wedge L_\alpha\models {\rm ZFC}\},$$
then by the same argument as above, $C$ is a required. The difference now is that if $\alpha\in C$ and $j\restr L_\alpha:L_\alpha\rightarrow L_{j(\alpha)}$, then $M\in {\rm Crit}(j)$ instead of just $\kappa={\rm crit}(j)$. \telos

\vskip 0.2in

The assumption that there is a  strongly critical model is no stronger  (over ZFC)  than  $Loc({\rm ZFC})$+ ``there is a measurable cardinal''.

\begin{Lem} \label{L:correspond}
{\rm ZFC+}$Loc({\rm ZFC})+${\rm ``$\kappa$ is a measurable cardinal''} proves that $V_\kappa$ is strongly critical. In particular the same is proved in {\rm ZFC+}{\rm ``there is a proper class of inaccessibles''+}  {\rm ``$\kappa$ is a measurable cardinal''}, as well as in the theory {\rm ZFC+}{\rm ``there is a proper class of $V_\alpha$ such that $V_\alpha\models{\rm ZFC}$+``$\kappa$ is a measurable cardinal''}.\footnote{In  \cite[Prop. 2.23]{Tz10} it is observed that  $Loc({\rm ZFC})$ is no stronger  than ``there is a proper class of  inaccessibles'', where the last principle is denoted by $IC^\infty$. Also in \cite[footnote 6]{Tz10} we denoted  by $NM$ the axiom ``there is a natural model'' (i.e., a $V_\alpha$ such that $V_\alpha\models{\rm ZFC}$). Let us denote by $NM^\infty$ the axiom ``there is a proper class of natural models''. It is  clear that over ZFC we have
$IC^\infty \Rightarrow NM^\infty \Rightarrow Loc({\rm ZFC})$.
Under mild large cardinal assumptions these arrows cannot be  reversed.}
\end{Lem}

{\em Proof.} Let $\kappa$ be measurable and let $x$ be a set. Let $j:V\rightarrow W$ be an elementary embedding of the universe with $\kappa={\rm crit}(j)$. By $Loc({\rm ZFC})$ there is a model $N\models{\rm ZFC}$ such that $\{V_\kappa,x\}\subset N$.  Let $j'=j\restr N$. Then $j(N)$ is also a model of ZFC and  $j':N\rightarrow j(N)$ is an elementary embedding with $V_\kappa\in {\rm Crit}(j)$. Thus $V_\kappa$ is $x$-strongly critical for every $x$. \telos

\vskip 0.2in

It seems  unlikely that one can  prove the inference of 4.10 without the assumption $Loc({\rm ZFC})$, i.e.,   to prove that strongly critical models exist in ZFC+``there is a measurable cardinal'' alone. This however can be done in  ZFC+``there is a strong  cardinal'' (see definition \ref{D:strongc}).

\begin{Lem} \label{L:alone}
{\rm (ZFC)} If $\kappa$ is a strong cardinal then $V_\kappa$ is strongly critical.
\end{Lem}

{\em Proof.} Let $\kappa$ be strong and let $x$ be a set. It suffices to show that $V_\kappa$ is $x$-strongly critical. Let $V_\alpha$ be such that $\{x,V_\kappa\}\subset V_\alpha$. Since $\kappa$ is $\alpha$-strong, there is a $j:V\rightarrow W$, where $W$ is an inner model, such that $j(\kappa)>\alpha$ and $V_{\kappa+\alpha}\subseteq W$. Then $j(\kappa)$ is inaccessible in $W$, therefore $W_{j(\kappa)}\models {\rm ZFC}$. Also for every $\beta\leq \kappa+\alpha$, $V_\beta=W_\beta$, therefore $V_\alpha=W_\alpha\subseteq W_{j(\kappa)}$. Let $K=W_{j(\kappa)}$. Then  $x,V_\kappa\in K$. If $j'=j\restr K$ and $N=W_{j^2(\kappa)}$, then the models $K$, $N$ and the elementary embedding $j':K\rightarrow N$ witness the fact that $V_\kappa$ is $x$-strongly critical.  \telos

\vskip 0.2in

The notions of critical and strongly critical model can be strengthened even further by demanding the elementary embedding to be an inner or self-embedding.

\begin{Def} \label{D:self}
{\em Let $M$ be a model and $x$ be a set.  $M$ is said to be} $x$-self-critical {\em (resp.} $x$-inner critical{\em)  if there is a  model $N$ and  an elementary embedding $j:N\rightarrow N$ (resp. if there are   models $N,K$ such that $K\subseteq N$ and  an elementary embedding $j:N\rightarrow K$), such that $\{M,x\}\subset N$ and  $M\in {\rm Crit}(j)$. $M$ is said to be} strongly self-critical {\em (resp.} strongly inner critical{\em) if it is $x$-self-critical (resp. $x$-inner critical) for every $x$.}
\end{Def}

The consistency strength of the existence of strongly inner-critical models is not  very high. Namely the following holds.

\begin{Lem}
{\rm (ZFC)} If $\kappa$ is a measurable cardinal and there exists a proper class of $V_\alpha$'s such that $V_\alpha\models{\rm ZFC}$ (i.e., if $NM^\infty$ holds), then $V_\kappa$  is a strongly inner critical model.
\end{Lem}

{\em Proof.} We have  to show that for every $x$ there is a model $M$ such that $\{V_\kappa, x\}\subseteq M$ and a $j:M\rightarrow N\subseteq M$ with $V_\kappa\in {\rm Crit}(j)$. Let  $U$ be a $\kappa$-complete ultrafilter on $\kappa$.  By $NM^\infty$ there is a $V_\alpha$ such that $\{V_\kappa,x,U\}\subseteq  V_\alpha$ and $V_\alpha\models{\rm ZFC}$. Then $\kappa$ is measurable in $V_\alpha$, so there is an elementary embedding $j:V_\alpha\rightarrow V_\alpha^\kappa/U\subseteq V_\alpha$ with ${\rm crit}(j)=\kappa$.  Therefore   $V_\kappa\in {\rm Crit}(j)$. \telos

\vskip 0.2in

However the consistency strength of the existence of strongly self-critical models seems to be much higher. We shall prove their consistency assuming the existence (in ZFC) of some kind of  rank-to-rank elementary embeddings. Recall that a cardinal $\kappa$ is  ${\rm I}_3$, denoted ${\rm I}_3(\kappa)$, if there is an elementary embedding $j:V_\lambda\rightarrow V_\lambda$ with ${\rm crit}(j)=\kappa$ (see \cite[p. 325]{Ka97}). ${\rm I}_3$ cardinals are $n$-huge for every $n<\omega$.  Yet we need something stronger.

\begin{Def} \label{D:rank}
{\em (ZFC) We say that a cardinal $\kappa$ is}  strongly ${\rm I}_3$, {\em denoted by ${\rm SI}_3(\kappa)$, if for every $\alpha>\kappa$ there is a $\lambda\geq \alpha$ and an elementary embedding $j:V_\lambda\rightarrow V_\lambda$ with  ${\rm crit}(j)=\kappa$.}
\end{Def}

\begin{Lem} \label{L:self}
{\rm (ZFC)} If  ${\rm SI}_3(\kappa)$, then $V_\kappa$ is strongly self-critical.
\end{Lem}

{\em Proof.} Suppose   ${\rm SI}_3(\kappa)$. Then, given  $x$, there is a  $\lambda> \kappa$, $rank(x)$, and an  elementary embedding $j:V_\lambda\rightarrow V_\lambda$ with  ${\rm crit}(j)=\kappa$. But then $V_\kappa\in {\rm Crit}(j)$. Hence $V_\kappa$ is $x$-self-critical for every $x$, thus $V_\kappa$ is strongly self-critical. \telos

\subsection{${\rm ZFC}+Loc({\rm ZFC})+\cdots$}
As already pointed out in \cite[p. 584]{Tz10}, ${\rm ZFC}+Loc({\rm ZFC})$ is a mild substitute of the theory ZFC + ``there is a proper class of inaccessible cardinals''. Accordingly,  ${\rm ZFC}+Loc^{\pi^1_1ind}({\rm ZFC})$ is a mild substitute of ZFC + ``there is a proper class of weakly compact  cardinals''. Mild here means ``local'', i.e., with no reference to large cardinals. More generally, for $\phi(x)$ a large  model property, the theories ${\rm ZFC}+Loc^\phi({\rm ZFC})$ and ${\rm ZFC}+ Loc({\rm ZFC})+\exists x \phi(x)$, or ${\rm ZFC}+ Loc({\rm ZFC})$ together with a  combination of the axioms $Loc^\phi({\rm ZFC})$ and $\exists x \phi(x)$, and  possibly augmented with ${\rm Found}_{On}$ or some piece of Replacement, seem to be worth  studying.
What would one expect to prove in such theories? The following could be among the expected results:

\vskip 0.1in

(a) Existence of  models of ZFC with special closure properties,  e.g. natural models.

(b) Existence of  ``internally strong'' models, i.e., satisfying ZFC+ a large cardinal property. (As it follows from  Lemma \ref{L:inac} (v), ${\rm ZFC}+ Loc({\rm ZFC})$+``there is an embeddable model'' proves that there is a model $M\models {\rm ZFC}$+``there is a strongly inaccessible cardinal''.)

(c) Existence (i.e., restoration) of large cardinals.

(d) Results having some impact on everyday mathematics. (Such results are least likely to be proved if  we judge by the analogous capabilities of the classical large cardinal axioms.)

\vskip 0.1in

The result below belongs to group (a) above.

\begin{Prop} \label{P:supercitical}
${\rm ZFC}+Loc({\rm ZFC})+${\rm ``there is a $V_\alpha$-strongly critical model $M$''}, for some $V_\alpha$ such that $M\subseteq V_\alpha$,    proves that {\em ``there is a natural model of ZFC''}.
\end{Prop}

{\em Proof.} Let $M$ be  $V_\alpha$-critical such that $M\subseteq V_\alpha$ in ZFC+$Loc({\rm ZFC})$, for some $\alpha$. Then  there are models $N,K$ and an elementary embedding $j:N\rightarrow K$ such that $M\in {\rm Crit}(j)$ and $V_\alpha\in N$. It follows that $N_\xi=V_\xi$ for every $\xi\leq \alpha$. Let $ht(M)=\beta$.  Since $M\subseteq V_\alpha$, we have that $\beta\leq\alpha$, therefore  $N_\beta=V_\beta$. Moreover ${\rm crit}(j)=\beta$ and  $N_\beta\in {\rm Crit}(j)$.  By Lemma \ref{L:inac} (v),   $N_\beta\models {\rm ZFC}$. It follows that $V_\beta\models {\rm ZFC}$, i.e.,   ZFC has a natural model. [Alternatively: Since $V_\beta\in {\rm Crit}(j)$, we have $V_\beta=N_\beta=K_\beta$. It follows that $j(V_\beta)=K_{j(\beta)}$ and  $V_\beta=K_\beta\in K_{j(\beta)}$, since $j(\beta)>\beta$. Thus $V_\beta\prec K_{j(\beta)}$ and $V_\beta\in K_{j(\beta)}$. So by Lemma \ref{L:montague}, $V_\beta\models{\rm ZFC}$.] \telos

\vskip 0.2in

The next result strengthens the previous one and belongs to group (c).

\begin{Prop} \label{P:toomuch}
${\rm ZFC+} Loc({\rm ZFC})+${\rm ``there is a strongly critical model''} proves that there is a strongly inaccessible cardinal.
\end{Prop}

{\em Proof.} Let $\kappa$ be a cardinal. By the Reflection principle of ZFC, the class of ordinals $\alpha$ such that for every $\kappa<\alpha$, $V_\alpha\models$ ``$\kappa$ is strongly inaccessible'' implies $\kappa$ is strongly inaccessible, is a proper class. So given a strongly critical model $M$, there is $\alpha$ such that $M\in V_\alpha$ and for every $\kappa$,
\begin{equation} \label{E:relative}
V_\alpha\models \mbox{``$\kappa$ is strongly inaccessible''} \ \Rightarrow \mbox{ $\kappa$ is strongly inaccessible}.
\end{equation}
Fix such a $V_\alpha$.  By strong criticalness of $M$, there are models  $N, K$ and an elementary embedding  $j:N\rightarrow K$ such that $V_\alpha\in N$ (so $\{M,V_\alpha\}\subset N$)  and  $M\in {\rm Crit}(j)$. Let $\kappa=ht(M)$. Then $\kappa={\rm crit}(j)$ and by Lemma \ref{L:inac} (v), $N\models$ ``$\kappa$ is strongly inaccessible''. Since $V_\alpha\in N$, a fortiori $V_\alpha\models$ ``$\kappa$ is strongly inaccessible''.
So by (\ref{E:relative}), $\kappa$ is strongly inaccessible. \telos

\vskip 0.2in

Let $scrit(x)$ be the  formula expressing ``$x$ is a strongly critical model''. A rather immediate consequence of the proof of \ref{P:toomuch} is the following.

\begin{Cor} \label{C:strong}
${\rm ZFC+}Loc^{scrit}({\rm ZFC})$ proves ``there is a proper class of inaccessible cardinals''.
\end{Cor}

\begin{Prop} \label{P:more}
${\rm ZFC+} Loc({\rm ZFC})+$ ``there is a strongly self-critical model'' proves that  there is a  cardinal $\kappa$ which is $n$-ineffable for every $n\in\omega$.
\end{Prop}

{\em Proof.} The proof is similar to that of \ref{P:toomuch}. Given a strongly self-critical $M$, let $\kappa=ht(M)$.  By Reflection there is $V_\alpha$ such that $M\in V_\alpha$,  $V_\alpha$ contains all mappings $f:[\kappa]^n\rightarrow \{0,1\}$, for every $n$, and
\begin{equation} \label{E:relative1}
V_\alpha\models ``(\forall n)(\mbox{$\kappa$ is $n$-ineffable)''} \ \Rightarrow (\forall n)(\mbox{ $\kappa$ is $n$-ineffable}).
\end{equation}
It is clear that if some function  $f:[\kappa]^n\rightarrow \{0,1\}$ has a stationary homogeneous  subset $X\subseteq \kappa$, then $X\in V_\alpha$.
Let $N$ be a model such that $V_\alpha\in N$. By strong self-criticalness there is an elementary embedding $j:N\rightarrow N$ such that $M\in {\rm Crit}(j)$. Since $\kappa=ht(M)$, it follows that  $\kappa\in {\rm crit}(j)$. By  Lemma  \ref{L:inac} (v), $N\models$ ``$\kappa$ is $n$-ineffable'' for every $n\in\omega$. Then, by the way we chose $N$ it is clear that $V_\alpha\models$ ``$\kappa$ is $n$-ineffable'' for every $n\in\omega$. So,  by (\ref{E:relative1}), $\kappa$ is $n$-ineffable for every $n\in\omega$. \telos

\vskip 0.2in

Let $sscrit(x)$ be the formula expressing ``$x$ is a strongly self-critical model''. As a consequence  of the proof of Proposition \ref{P:more} we have:

\begin{Cor} \label{C:strong1}
${\rm ZFC+}Loc^{sscrit}({\rm ZFC})$ proves ``there is a proper class of $n$-ineffable cardinals'', for every $n\in\omega$.
\end{Cor}

\section{LZFC and  $V=L$. The Tall Model Axiom}
In this section we examine the relationship of LZFC with $V=L$ and especially with $V\neq L$.
We shall need below  the following well-known facts of ZFC (see e.g. \cite[Th.18.20]{Je03}, \cite[Th. 4.3]{De84} and \cite[Th. 9.12, 9.17]{Ka97}).

\begin{Thm} \label{T:Devlin}
{\rm (ZFC)} (i) The following are equivalent:

(a) $0^\#$ exists.

(b) There is an elementary embedding $j:L\rightarrow L$.

(c) There is an elementary embedding $j:L_\alpha\rightarrow L_\beta$, where $\alpha,\beta$ are limit ordinals, with ${\rm crit}(j)=\kappa<|\alpha|$.

(ii) If $0^\#$ exists, then there is a proper class of ordinals $\alpha$ such that $L_\alpha\prec L$.

(iii) If  $0^\#$ exists, then for every infinite $x\in L$, $|{\cal P}(x)^L|=|x|$.
\end{Thm}

First concerning the consistency of ${\rm LZFC}+V=L$, let us observe the following:

\begin{Lem} \label{L:atfirst}
(i) If {\rm ZFC+``there is an inaccessible''} is consistent, then so is {\rm ZFC+}$Loc({\rm ZFC})+V=L$.

(ii) If {\rm ZFC+``$0^\#$ exists''} is consistent, then so is {\rm ZFC+}$Loc^{sext}({\rm ZFC})+V=L$.
\end{Lem}

{\em Proof.} (i) It was seen in \cite[Pr. 2.23 (ii)]{Tz10} that if $\kappa$ is a strongly inaccessible cardinal in the ZFC universe $V$, then $V_\kappa\models {\rm ZFC}+Loc({\rm ZFC})$. From this it is straightforward that $L_\kappa\models{\rm ZFC}+Loc({\rm ZFC})+V=L$.

(ii) Let $M$ be a model of  ZFC+``$0^\#$ exists''. It suffices to see that  $L^M\models Loc^{sext}({\rm ZFC})$. This follows from the fact that, by \ref{T:Devlin} (ii),  there is in $M$ a cofinal class of models $L_\alpha$ such that $L_\alpha\prec L^M$. Then,  as in the proof of Proposition \ref{P:nostrong}, all these models are strongly extendible, so $L^M\models Loc^{sext}({\rm ZFC})$.   \telos

\vskip 0.2in

In view of the fact that existence of $0^\#$ is the standard means to refute $V=L$ over ZFC, the question is what  ``large model axioms'' are needed  in order to refute $V=L$ over LZFC.  By \ref{L:atfirst} (ii), $V=L$ is consistent even with strengthenings of $Loc({\rm ZFC})$, like $Loc^{sext}({\rm ZFC})$.

So below we shall try  some new axioms  giving information about  the internal truths of  local models, e.g., about how they see the cardinalities of certain sets. For instance suppose a set $x$ is given. It is a straightforward consequence of $Loc({\rm ZFC})$ that  there is a set $y$ and  a model $M$ of ZFC such that $\{x,y\}\subset  M$ and $M\models |x|<|y|$. For simplicity, and without loss of generality we can deal just  with ordinals  instead of arbitrary sets. Thus  the statement $$(\forall \kappa)(\exists \alpha>\kappa)(\exists M)(\alpha\in M \wedge  M\models |\kappa|<|\alpha|)$$
is a theorem of LZFC. Now let us make the question a bit harder. Let $\kappa$ be given. Does there exist an $\alpha>\kappa$ such that for every $\delta\geq  \alpha$ there is a model $M$ such that $\delta\in M$ (equivalently $ht(M)>\delta$) and $M\models |\alpha|>|\kappa|$?
This is the statement:
\begin{equation} \label{E:presaxiom}
(\forall \kappa)(\exists \alpha>\kappa)(\forall \delta\geq \alpha)(\exists M)(\delta\in M \wedge  M\models |\kappa|<|\alpha|).
\end{equation}
(\ref{E:presaxiom}) says that for every $\kappa$ there is an $\alpha>\kappa$ such that there are arbitrarily ``tall'' models which do not collapse $\alpha$ to $\kappa$. We shall refer to (\ref{E:presaxiom}) as the {\em Tall Model Axiom}, or $TMA$ for short.
$TMA$ is a  principle of mild consistency strength as is seen by the following:

\begin{Lem} \label{L:cons}
(i)  ${\rm ZFC}+Loc({\rm ZFC})\vdash TMA$.

(ii) If $\lambda$ is a limit cardinal in ${\rm ZFC}+Loc({\rm ZFC})$, then $H(\lambda)\models {\rm LZFC}+TMA$. More generally, if $N\models {\rm LZFC}$ and $N$ does not have a greatest cardinality, then $N\models TMA$.
\end{Lem}

{\em Proof.} (i)  Given $\kappa$, just pick an (absolute)  cardinal $\alpha$ such that $\kappa<\alpha$. For every $\delta\geq \alpha$  pick, by $Loc({\rm ZFC})$, a  model $M$ such that $\delta\in M$. Then obviously $M\models |\kappa|<|\alpha|$.

(ii) Let $\lambda$ be a limit cardinal in the theory ${\rm ZFC}+Loc({\rm ZFC})$.  That $H(\lambda)$ satisfies $Loc({\rm ZFC})$  follows from $Loc({\rm ZFC})$. Indeed, let $x\in H(\lambda)$, i.e.,  $|TC(x)|<\lambda$. By $Loc({\rm ZFC})$ there is a  model $M$ of ZFC such that $x\in M$. Then, by L\"{o}wenheim-Skolem there is an elementary submodel $N$ of $M$ containing  $TC(x)$ and $|N|<\lambda$. If $\bar{N}$ is the  transitive collapse of $N$, then $x\in \bar{N}$, ${\bar N}\models {\rm ZFC}$ and   $\bar{N}\in H(\lambda)$. Now let $\kappa_0\in H(\lambda)$. Since $\lambda$ is a limit cardinal, there is a cardinal $\mu\in H(\lambda)$ such that $\kappa_0<\mu$. For every ordinal  $\delta\geq \mu$ pick as before, by $Loc({\rm ZFC})$, a model $K$ such that $\delta\in K$. Then obviously $K\models |\kappa_0|<|\mu|$. By L\"{o}wenheim-Skolem $K$ can be chosen  so that  $|K|<\lambda$. Therefore $H(\lambda)\models TMA$. The proof of the more general statement is the same. \telos

\vskip 0.2in

We can now  see that if  $TMA$ is added to the theory LZFC+``there exists a strongly critical model'', then $V=L$ fails.

\begin{Thm} \label{P:holds}
{\rm LZFC+}TMA{\rm +``there exists a strongly critical model''} proves   $V\neq L$.
\end{Thm}

{\em Proof.}  We assume {\rm LZFC}+$TMA$+``there exists a strongly critical model''+$V=L$. It suffices to reach a contradiction. Let $M$ be strongly critical with $ht(M)=\kappa_0$. By $TMA$ there is an $\alpha_0>\kappa_0$ such that
\begin{equation} \label{E:presaxiom1}
(\forall \delta\geq \alpha_0)(\exists N)(\delta\in N \wedge  N\models |\kappa_0|<|\alpha_0|).
\end{equation}
Since $V=L$, $M\in L$, so by proposition \ref{P:Lalpha} (ii)  there are  limit ordinals  $\beta_0, \gamma_0$ such that $\alpha_0<\beta_0\leq \gamma_0$ and an elementary embedding $j:L_{\beta_0}\rightarrow L_{\gamma_0}$ with $M\in {\rm Crit}(j)$. Hence ${\rm crit}(j)=\kappa_0$. By $V=L$, $j\in L$. Therefore $j\in L_\delta$ for some $\delta>\gamma_0$. Then  it follows  from (\ref{E:presaxiom1}) that there is a model $N$ of ZFC which contains $\delta$ and  $N\models |\kappa_0|<|\alpha_0|$, so  $N\models |\kappa_0|<|\beta_0|$. A fortiori, $L^N\models |\kappa_0|<|\beta_0|$. Moreover  $j\in L^N$ since $j\in L_\delta\in L^N$. Thus   $L^N$ is a model of ZFC that contains an elementary embedding $j:L_{\beta_0}\rightarrow L_{\gamma_0}$ with critical point $\kappa_0$ such that $L^N\models |\kappa_0|<|\beta_0|$. But then also $L^N\models \kappa_0<|\beta_0|$, since the ordinals $\kappa_0$ and $|\kappa_0|^{L^N}$ are equipollent in $L^N$. By Theorem \ref{T:Devlin} (i), $L^N\models ``0^\# \ \mbox{exists''}$. This contradicts the fact that  $L^N\models V=L$. \telos

\begin{Quest} \label{Q:remove}
Can we remove $TMA$ from the assumptions of the previous theorem, i.e., does  {\rm LZFC+``there exists a strongly critical model''} prove   $V\neq L$?
\end{Quest}

We have already seen in Lemma \ref{L:cons} (ii)  that $TMA$ holds in any  universe of LZFC without greatest cardinality. So let us focus on this last property.  Recall that an ordinal $\alpha$  is said to be a cardinal in LZFC, if (exactly as in ZFC) there is no $\beta<\alpha$ such that $\alpha\sim \beta$. Moreover, the (absolute) cardinal of $x$, denoted $|x|$, is the least $\alpha$ (if there exists) such that $x\sim \alpha$. However, due to the lack of $\in$-induction, we cannot ensure that for any  given set $x$  $|x|$ exists.   But the  symbol $|x|$ can keep being  used  in the expressions $|x|=|y|$, $|x|\leq |y|$,  $|x|<|y|$ which have the ordinary meaning (that is, existence of bijections or injections between $x,y$, etc.) Below we shall denote by  $GC$ the statement ``there is a greatest cardinality'', i.e.,
$$(GC) \hspace{.5\columnwidth minus .5\columnwidth}(\exists \kappa)(\forall y)(|y|\leq |\kappa|).\hspace{.5\columnwidth minus .5\columnwidth} \llap{}$$
A particular case of $GC$ is the axiom $(\forall x)(|x|\leq \omega)$, which we often refer to as ``all sets are countable''.

\begin{Lem} \label{L:nodif}
The following are equivalent over {\rm LZFC}:

(i) $GC$ (:= $(\exists \kappa)(\forall y)(|y|\leq |\kappa|))$.

(ii) $(\exists \kappa)(\forall \alpha)(|\alpha|\leq |\kappa|)$.

(iii) $(\exists x)(\forall \alpha)(|\alpha|\leq |x|)$.

(iv) $(\exists x)(\forall y)(|y|\leq |x|)$.
\end{Lem}

{\em Proof.} (i)$\Rightarrow$ (ii) and  (ii)$\Rightarrow$ (iii)  are trivial.

(iii)$\Rightarrow$(iv): Assume $x_0$ is a set satisfying (iii) and let $y$ be any  given set. We have to show that $|y|\leq |x_0|$.  Working in a model $M$ containing $y$, we can find an ordinal $\alpha$ such that  $|\alpha|=|y|$. Then  by (iii) $|\alpha|\leq |x_0|$. Therefore $|y|\leq |x_0|$.

(iv)$\Rightarrow$(i): Assume  $x_0$ is a set satisfying (iv) and find an ordinal $\kappa$ as before such that $|x_0|=|\kappa|$. Let $y$ be any given set. By (iv) $|y|\leq |x_0|$, therefore $|y|\leq |\kappa|$.
 \telos

\begin{Lem} \label{L:greatcard}
${\rm LZFC}+\neg(GC)\vdash TMA$.
\end{Lem}

{\em Proof.} This has essentially been shown in Lemma \ref{L:cons} (ii). Suppose in LZFC there is no greatest cardinal, and let $\kappa$ be an ordinal. Then there is an ordinal $\alpha>\kappa$ such that $|\kappa|<|\alpha|$. Then for every $\delta\geq \alpha$ and every model $M$ such that $\delta\in M$, obviously $M\models|\kappa|<|\alpha|$. Therefore $TMA$ is true. \telos

\vskip 0.2in

The converse of \ref{L:greatcard} is open.

\begin{Quest} \label{Q:consistent}
Is   ${\rm LZFC}+ GC+TMA$  consistent?
\end{Quest}

What we can only prove here is the following:

\begin{Prop} \label{P:allcountable}
${\rm LZFC}+GC+V=L\vdash \neg(TMA)$. \\
In particular,

${\rm LZFC}+ \ \mbox{\rm ``every set is  countable''}+V=L\vdash \neg(TMA)$.
\end{Prop}

{\em Proof.} Assume  ${\rm LZFC}+GC+V=L$. We have to show $\neg(TMA)$, i.e., that
\begin{equation} \label{E:negation}
(\exists \kappa)(\forall  \alpha>\kappa)(\exists \delta)(\forall M)(\delta\in M
\Rightarrow M\models|\kappa|=|\alpha|)
\end{equation}
is true. By assumption there is an ordinal  $\kappa$ such that $\forall \alpha (|\alpha|\leq |\kappa|)$. Fix such a $\kappa_0$. It suffices to  show that
\begin{equation} \label{E:needed}
(\forall \alpha>\kappa_0)(\exists \delta)(\forall M)(\delta\in M\Rightarrow M\models |\kappa_0|=|\alpha|).
\end{equation}
Pick an  $\alpha>\kappa_0$. Since $|\alpha|\leq |\kappa_0|$,  there is (in $V$) a bijection $f:\kappa_0\rightarrow \alpha$. Since $V=L$, $f\in L$. Hence there is a $\delta$ such that  $f\in L_\delta$. Then  $\delta$ witnesses the truth of (\ref{E:needed}), since if  $M$ is a model of ZFC such that $\delta\in M$, then $L_\delta\in M$, and so $f\in M$. Therefore $M\models |\alpha|=|\kappa_0|$. \telos

\vskip 0.1in

It follows immediately from \ref{L:greatcard} and \ref{P:allcountable} that $\neg(TMA)$ and $GC$ are equivalent over LZFC+$V=L$, i.e.,

\begin{Cor} \label{C:overL}
${\rm LZFC}+V=L\vdash \neg(TMA)\Leftrightarrow GC$.
\end{Cor}

\begin{Cor} \label{C:cons}
(i)  ${\rm LZFC}+ \neg(GC) + \neg(TMA)$ is inconsistent.

(ii)  ${\rm LZFC}+ \neg(GC) + TMA$ is consistent, provided  ${\rm ZFC}+Loc({\rm ZFC})$ is so.

(iii) ${\rm LZFC}+ GC+\neg(TMA)$ is consistent, provided ${\rm ZFC}+Loc({\rm ZFC})+V=L$ is so.
\end{Cor}

{\em Proof.} (i) follows from Lemma \ref{L:greatcard}. (ii) follows from Lemma \ref{L:cons} (ii). (iii) follows from Proposition  \ref{P:allcountable}. \telos

\vskip 0.2in

Before closing this section  let us consider the following  reasonable variants of  $TMA$:

$$(TMA_1) \hspace{.5\columnwidth minus .5\columnwidth}(\forall x)(\exists \alpha)(\forall \delta\geq \alpha)(\exists M)(\{x,\delta\}\subset M \ \wedge \  M\models |x|<|\alpha|).\hspace{.5\columnwidth minus .5\columnwidth} \llap{}$$

$$(TMA_2) \hspace{.5\columnwidth minus .5\columnwidth}(\forall x)(\exists y)(\forall \delta)(\exists M)(\{x,y,\delta\}\subset M \ \wedge \ M\models |x|<|y|).\hspace{.5\columnwidth minus .5\columnwidth} \llap{}$$

Concerning the relative   strength of $TMA$, $TMA_1$ and $TMA_2$, it is immediate that $TMA_1\Rightarrow TMA$ and $TMA_1\Rightarrow TMA_2$, over LZFC.

\begin{Lem} \label{L:greatcard1}
${\rm LZFC}+\neg(GC)\vdash TMA_1$, therefore ${\rm LZFC}+\neg(GC)\vdash TMA_2$.
\end{Lem}

{\em Proof.} Similar to that of Lemma \ref{L:greatcard}. \telos

\begin{Prop} \label{P:greatcard1}
(i) ${\rm LZFC}+GC+V=L\vdash \neg(TMA_2)$.

(ii) ${\rm LZFC}+GC+V=L\vdash \neg(TMA_1)$.
\end{Prop}

{\em Proof.} Similar to that of Proposition \ref{P:allcountable}. \telos

\begin{Cor} \label{L:alleq}
$TMA$, $TMA_1$, $TMA_2$ and $\neg(GC)$ are all equivalent over ${\rm LZFC}+V=L$.
\end{Cor}

{\em Proof.} It follows immediately from \ref{L:greatcard1} and \ref{P:greatcard1}. \telos

\vskip 0.2in

Finally we include  here a result of independent interest and potential applicability, saying that $Loc({\rm ZFC})$ is preserved under forcing extensions.\footnote{Generic extensions preserve also the related principles $NM^\infty$ and $IC^\infty$ mentioned in footnote 5 above.}

\begin{Prop} \label{P:include}
Any  generic extension of a model of ${\rm ZFC}+Loc({\rm ZFC})$ is also a model  of $Loc({\rm ZFC})$.
\end{Prop}

{\em Proof.} Let $M\models{\rm ZFC}+Loc({\rm ZFC})$, $\P\in M$ be a forcing notion, and $G\subseteq \P$ be $M$-generic. To show that $M[G]\models Loc({\rm ZFC})$, pick $a\in M[G]$. It suffices to find  $M_1\in M[G]$ such that $M_1\models {\rm ZFC}$ and $a\in M_1$. Let $\dot{a}\in M$ be a $\P$-name of $a$. Since $M\models Loc({\rm ZFC})$ there is  $N\in M$ such that $\{\P, \dot{a}\}\subset N$ and $N$ is a  model of ZFC in the sense of $M$. Now $N\subseteq M$ and $G$ is $M$-generic, therefore    $G$  is also $N$-generic, so   $N[G]\models{\rm ZFC}$. Moreover,  $\dot{a}$ is a $\P$-name in the sense of $N$, so $a\in N[G]$ and $N[G]\in M[G]$. Thus $N[G]$ is the required model. \telos

\section{Vop\v{e}nka's Principle}
Among strong large cardinal properties Vop\v{e}nka's Principle ($VP$)  seems to be particularly fitting  to the context of LZFC. In this section we shall show that if we add $VP$ to a strengthened version of LZFC, then ZFC is recovered.

Recall that $VP$ is a scheme rather than a single axiom, defined as follows: Given a formula $\phi(x)$ in one free variable, let $X_\phi$ denote the extension $\{x:\phi(x)\}$ of $\phi$. Then clearly ``$X_\phi$ is a proper class'' is a first-order statement. In particular the following expression is a first-order statement:

\vskip 0.1in

$(VP_\phi$) \ \ If $X_\phi$ is a proper class  of structures (of some fixed \\  \hspace*{0.8in} first-order language),  then there  are distinct $x,y\in X_\phi$
\\ \hspace*{0.8in} and an elementary  embedding $j:x\rightarrow y$,

\vskip 0.1in

\noindent (where  $j:x\rightarrow y$ may be trivial, i.e., $x\prec y$). Then  $VP=\{VP_\phi:\phi\in {\cal L}\}$. That $VP$ is an appealing scheme for LZFC follows from the fact that, as a consequence of $Loc({\rm ZFC})$, every set belongs to  a proper class of models of ZFC. Indeed, for every set $x$, let ${\cal M}(x)$  be the class of models of ZFC  containing $x$. That is,
$${\cal M}(x)=\{y: x\in y \ \wedge Tr(y) \ \wedge \ (y,\in)\models {\rm ZFC}\}.$$

\begin{Lem} \label{L:properclass}
{\rm (LZFC)}  For every $x$, ${\cal M}(x)$ is a  proper class.
\end{Lem}

{\em Proof.}  $Loc({\rm ZFC})$ implies that for every $x$,  $\bigcup{\cal M}(x)=V$.  Thus if ${\cal M}(x)$ is a set,  so is $V$. But then, by $\Delta_0$-Separation, so is also $R=\{x\in V:x\notin x\}$, that is, Russell's paradox reappears. \telos

\vskip 0.2in

The consequences  of $VP$ established  below  stem from two sources:

Firstly,  the  fact that in a proper-class context $VP$ works  as a  set existence principle. Specifically, $VP$ is an implication of the form ``if $X_\phi$ is a proper class, then such and such is the case''. Taking the contrapositive we have equivalently, ``if  such and such is not the case, then $X_\phi$ is a set''.

Secondly, an old theorem of  P. Vop\v{e}nka, A. Pultr and Z. Hedrl\'{i}n concerning  rigid graphs. Given a set $A$ and a binary relation $R\subset A\times A$, we refer to $(A,R)$ as a {\em graph.} Given a graph $(A,R)$, a mapping $f:A\rightarrow A$ is  an {\em endomorphism} if for all $x,y\in A$, $(x,y)\in R\Rightarrow (f(x),f(y))\in R$. The graph $(A,R)$ is said to be {\em rigid}  if the only endomorphism of $(A,R)$ is the identity.   The  Vop\v{e}nka-Pultr-Hedrl\'{i}n (V-P-H) result \cite{VPH65} is the following:

\begin{Thm} \label{T:keyVopenka}
{\rm (V-P-H \cite{VPH65})} {\rm (ZFC)}
For any infinite set  $A$, there is a binary  relation $R\subset A\times A$ such that the graph $(A,R)$ is rigid.\footnote{Theorem \ref {T:keyVopenka}  provides a proper class whose members are {\em individually} rigid structures. As  reported  in \cite{AR94}, pp. 278-279, after the proof of this result  P. Vop\v{e}nka believed
that not only such a proper class of individually rigid structures, but also a proper class of structures which is {\em rigid in itself}  could  be constructed. That is to say, he believed that there exists  a proper class of structures such that no non-trivial homomorphism exists between any two distinct elements of it  (which is, roughly, the negation of what we call today $VP$). He then playfully set forth  the {\em negation} of this conjecture, i.e., what we now call $VP$, just  in order to ``tease'' other mathematicians, and make them finally disprove something he believed  was definitely false.} \label{F:4}
\end{Thm}
Henceforth we refer to Theorem \ref{T:keyVopenka} as  ``V-P-H''.
V-P-H was originally proved in ZFC. It is open whether  it is provable in LZFC, but as we shall see there  are reasonable ``local style'' extensions of LZFC that prove it.  The next theorem is the main result of this section.

\begin{Thm} \label{T:main}
If   $T$ is a theory such that ${\rm LZFC}\subseteq T$ and $T\vdash $V-P-H, then ${\rm ZFC}\subseteq T+VP$.
\end{Thm}

{\em Proof.} Let $T$ be as stated. It suffices to show that $T+VP$ proves Replacement and Powerset. So the theorem is an immediate consequence of the next two Lemmas. \telos

\begin{Lem} \label{L:Delta}
Let  $T$ be a theory such that ${\rm LZFC}\subseteq T$ and $T$ proves V-P-H. Then:

(i) $T+VP\vdash \Delta_0$-{\rm Repl}.

(ii) $T+VP+\Delta_0${\rm -Repl}$ \ \vdash$ {\rm Repl}.

(iii) $T+VP\vdash$ {\rm Repl}.
\end{Lem}

{\rm Proof.} (i) Let  $T$ be as stated and let us  work in $T+VP$.  To prove $\Delta_0$-Replacement, let $\phi(x,y)$ be a $\Delta_0$ formula such that $(\forall x)(\exists !y)\phi(x,y)$. As usual this  defines a class mapping $F_\phi:V\rightarrow V$ such that $F_\phi(x)=y$ iff $\phi(x,y)$. Fix a set $A$. We have to show that the class $B=F_\phi''A=\{F_\phi(x):x\in A\}$ is a set.  By assumption V-P-H is true in $T$, so fix a binary  relation $R\subset A\times A$ such that $(A,R)$ is rigid. For every $b\in B$ consider the structure $K_b=(A\times \{b\}, R_b, E_b)$, where $R_b$ is the binary relation on $A\times \{b\}$ induced by $R$, i.e., $$((x,b),(y,b))\in R_b\iff (x,y)\in R,$$  and $E_b$ is the unary relation on $A\times \{b\}$ defined by:
$$(x,b)\in E_b\iff F_\phi(x)=b.$$ Let ${\cal Z}=\{K_b:b\in B\}$. We claim that it suffices to show that ${\cal Z}$ is a set. For suppose  ${\cal Z}$ is a set. Then clearly  for some $n\in\N$, $B\subseteq \bigcup^n{\cal Z}=s$ and $s$ is a set. Therefore,
$$B=\{y\in s:(\exists x\in A)(F_\phi(x)=y)\}=\{y\in s:(\exists x\in A)\phi(x,y)\}.$$
Since the formula $(\exists x\in A)\phi(x,y)$ is $\Delta_0$ and $\Delta_0$-Separation holds in $T$ (because ${\rm LZFC}\subseteq T$ and $\Delta_0$-Separation holds in LZFC),  it follows that $B$ is  set.

So let us verify that ${\cal Z}$ is a set. Towards reaching a contradiction assume that ${\cal Z}$ is a proper class.  Then by $VP$ there are $b,c\in B$, $b\neq c$,  and an elementary embedding $f:K_b\rightarrow K_c$. $f$ induces the  mapping  $g:A\rightarrow A$ such that for every $x\in A$, $f(x,b)=(g(x),c)$. It is easy to see that $g:(A,R)\rightarrow (A,R)$ is an endomorphism. Indeed, by the elementarity of $f$ we have,
$$(x,y)\in R\Leftrightarrow ((x,b),(y,b))\in R_b\Leftrightarrow (f(x,b),f(y,b))\in R_c\Leftrightarrow$$$$ ((g(x),c),(g(y),c))\in R_c\Leftrightarrow (g(x),g(y))\in R.$$
By the rigidity of $(A,R)$, $g=id_A$, so for every $x\in A$,  $f(x,b)=(x,c)$. But then, for every $x\in A$ we have:
$$F_\phi(x)=b\Leftrightarrow (x,b)\in E_b\Leftrightarrow f(x,b)\in E_c\Leftrightarrow (x,c)\in E_c\Leftrightarrow F_\phi(x)=c,$$
a contradiction, since $b\neq c$.

(ii) Now we work in $T+VP$+$\Delta_0$-Repl, and prove that full Replacement holds. The proof is to a great extent similar to that of clause (i) above. Let $\phi(x,y)$ be a formula such that $(\forall x)(\exists !y)\phi(x,y)$, and let $F_\phi(x)=y$ iff $\phi(x,y)$. We fix again a set $A$ and show that if $B=F_\phi''A$, then $B$ is a set. We define the structures $K_b$ as before and set  ${\cal Z}=\{K_b:b\in B\}$. The only departure from the proof of (i) is on how we infer that $B$ is a set from ${\cal Z}$ being a set. Since $\phi$ now need not be $\Delta_0$ and full Separation is not available in LZFC,  the previous argument does not work. But we can invoke $\Delta_0$-Replacement. In view of the latter,  in order to show that $B$ is a set  it suffices to show that ${\cal Z}$ is a set, since  the  mapping ${\cal Z}\ni K_b\mapsto b\in B$ is  $\Delta_0$-definable and onto. Then we continue exactly as in (i). That is, we  assume that ${\cal Z}$ is a proper class and  reach the same contradiction.

(iii) Immediate from (i) and (ii). \telos

\begin{Lem} \label{L:rigid}
Let  $T$ be a theory such that ${\rm LZFC}\subseteq T$ and $T$ proves V-P-H. Then
$T+VP+\Delta_0${\rm -Repl}$ \ \vdash {\rm Powerset}$. Therefore, by \ref{L:Delta} (i), $T+VP\vdash {\rm Powerset}$.
\end{Lem}

{\em Proof.}  Let $T$ be as stated and assume the contrary, i.e., in  $T+VP+\Delta_0$-Repl  there is an infinite set $A$ such that ${\cal P}(A)$ is a proper class. It suffices to reach a contradiction.  Fix such a set $A$.  Fix also, by theorem V-P-H, a binary relation $R\subset A\times A$ such that $(A,R)$ is rigid. For every set $X\in {\cal P}(A)$ consider the first-order structure $S_X=(A,R,X)$, where $X$ is interpreted as a unary predicate. Let ${\cal X}=\{S_X:X\in {\cal P}(A)\}$.
The mapping ${\cal X}\ni S_X\mapsto X\in {\cal P}(A)$ is $\Delta_0$-definable and onto, so since ${\cal P}(A)$ is a proper class,  ${\cal X}$ is a proper class too.  By $VP$, there are $X,Y\in {\cal P}(A)$,  $X\neq Y$,  and an elementary embedding $f:S_X\rightarrow S_Y$, i.e., $f:(A,R,X)\rightarrow (A,R,Y)$. Then in particular  $f:(A,R)\rightarrow (A,R)$ is an endomorphism, so $f=id_A$ since $(A,R)$ is rigid. Now  for every $x\in A$, $S_X\models X(x) \Leftrightarrow S_Y\models Y(f(x))$, or $x\in X\Leftrightarrow f(x)\in Y$. Since $f=id_A$, it follows  that $X=Y$, a contradiction. \telos

\vskip 0.2in

{\em Proof of Theorem \ref{T:main}}. In view of Lemmas \ref{L:Delta} and \ref{L:rigid}, the proof of Theorem \ref{T:main} is complete. \telos

\begin{Rem} \label{R:Pomega}
{\em Note that if instead of the general Powerset axiom   one wants only to show that  ${\cal P}(\omega)$ is a set, one  can prove this simply in ${\rm LZFC}+VP+\Delta_0${\rm -Repl}, skipping Theorem V-P-H  and using just the rigidity of $\omega$. Namely, it suffices to consider for every  $X\subseteq \omega$, the structure $K_X=(\omega,S,0,X)$, where $S$ is the successor operation of $\omega$ and $X$ is a unary predicate. Then working exactly as in the proof of \ref{L:rigid}, and using the fact that $(\omega,S,0)$ is rigid, one  concludes that  ${\cal P}(\omega)$ is a set. }
\end{Rem}

Let us next turn to the question: Is LZFC sufficient for the proof of  V-P-H? This requires an inspection of the ZFC proof (see also \cite[Lemma 2.64]{AR94} for  a more concise version of it). The  idea of the proof is, given a set $A$, first  to identify  $A$ with an ordinal of the form $\lambda+2$  (using   $AC$), and then to define  $R$ by cleverly   employing  the ordinals $\alpha<\lambda+2$ with {\em countable cofinality.}
Then the rigidity of  $(\lambda+2,R)$, i.e., the fact that for every endomorphism $f$, $f(\beta)=\beta$, for every $\beta<\lambda+2$,  is shown by induction on $\lambda+2$. Thus the property ``${\rm cf}(\alpha)=\omega$'' plays a crucial role in the construction of $R$. This property is expressed by a $\Sigma_1$ sentence, so the construction is not absolute. It means that if we work in  LZFC and given a set $A$, we pick a model $M$ such that $A\in M$ and define  a relation $R$ on $A$ such that $M\models$ ``$(A,R)$ is rigid'', we cannot be sure  that $(A,R)$ is actually rigid. So LZFC does not seem adequate to prove V-P-H, at least {\em in the way} the latter is proved in ZFC.\footnote{Of course the possibility of a {\em new} proof of V-P-H in LZFC,  essentially  different from that of  \cite{VPH65}, cannot be excluded.}  But if given a set $A$, we can pick a model $M$ such that $A\in M$ and $M\models {\rm cf}(\alpha)=\omega$ iff  ${\rm cf}(\alpha)=\omega$, then the construction of $R$ in $M$ guarantees the  absolute rigidity of $(A,R)$.

Thus we arrive at the following definition which looks rather {\em ad hoc} but gives precisely the minimum requirements for an extension of LZFC in order to prove  V-P-H.

\begin{Def} \label{D:special}
{\em (LZFC) A model $M$ is said to be} special {\em if for every ordinal $\alpha\in M$, ${\rm cf}^M(\alpha)=\omega\Leftrightarrow {\rm cf}(\alpha)=\omega$.\footnote{This particular kind of ``special models'' is not connected with special models as used in classical model theory.}}
\end{Def}
Let $spec(x)$ be the formal expression of the property ``$x$ is a special model''. Let as usual ${\rm LZFC}^{spec}={\rm LZFC}+Loc^{spec}({\rm ZFC})$.

\begin{Thm} \label{T:special}
${\rm LZFC}^{spec}\vdash$ V-P-H. Therefore ${\rm ZFC}\subseteq {\rm LZFC}^{spec}+VP$.
\end{Thm}

{\em Proof.} Let $A$ be a set. By $Loc^{spec}({\rm ZFC})$ we can pick a  special model $M$ such that $A\in M$. Applying V-P-H inside $M$ we construct  an $R\subset A\times A$ as in the standard proof so that $(A,R)$  is rigid in $M$. Namely we identify $A$ with an ordinal $\lambda+2$ of $M$ and define $R$ as in V-P-H.   We claim that $(\lambda+2,R)$ is absolutely rigid. Indeed, let $f$ be an endomorphism of $(\lambda+2,R)$ outside $M$. By $Loc({\rm ZFC})$ there is a model $N$ such that $M\subset N$ and $f\in N$. Since $M$ is special, for every $\alpha<\lambda+2$, $N\models{\rm cf}(\alpha)=\omega$ iff $M\models{\rm cf}(\alpha)=\omega$ iff ${\rm cf}(\alpha)=\omega$. Thus in $N$ one can show by induction that  $f(\beta)=\beta$ for every $\beta<\lambda+2$,  exactly as this is proved in $M$ for the endomorphisms of $M$. It follows that $f=id$.
The other claim follows from Theorem \ref{T:main}. \telos

\vskip 0.2in

Now special models are related to strongly extendible models considered in section 3, and in particular their property to be $\Sigma_1$-elementary submodels of $V$ (see Lemma \ref{L:abssigma} (i)). As already remarked  in footnote \ref{f:4}, not the full strength of strong extendibility is needed for that property but only  the $\Sigma_1$ part  of it. So the following is an intermediate property between being special and strongly extendible.

\begin{Def} \label{D:sigmaext}
{\em A model $M$ is said to be}  $\Sigma_1$-strongly extendible {\em  if for every $x$ there is a model $N$ such that $x\in N$ and $M\prec_{\Sigma_1}N$.}
\end{Def}

Let  $\Sigma_1sext(x)$ denote  the formula expressing  ``$x$ is $\Sigma_1$-strongly extendible''.

\begin{Lem} \label{L:specsigma}
(i) If $M$ is $\Sigma_1$-strongly extendible, then $M\prec_{\Sigma_1}V$.

(ii) For every $M$,
 $$sext(M)\Rightarrow\Sigma_1sext(M)\Rightarrow spec(M).$$

(iii) ${\rm LZFC}^{spec}\subseteq  {\rm LZFC}^{\Sigma_1sext}\subseteq {\rm LZFC}^{sext}$.
\end{Lem}

{\em Proof.} (i) Let $M$ be $\Sigma_1$-strongly extendible, and let $\exists \overline{x}\phi(\overline{x},\overline{c})$ be a $\Sigma_1$ sentence with $\overline{c}\in M$ such that $V\models \exists \overline{x}\phi(\overline{x},\overline{c})$. Then $V\models \phi(\overline{a},\overline{c})$ for some $\overline{a}$. By $\Sigma_1$-strong extendibility there is a model $N$ such that $\overline{a},\overline{c}\in N$, and $M\prec_{\Sigma_1} N$. Then  $N\models \phi(\overline{a},\overline{c})$ since $\phi$ is bounded. Thus $N\models \exists \overline{x}\phi(\overline{x},\overline{c})$, so  $M\models \exists \overline{x}\phi(\overline{x},\overline{c})$ too.

(ii) The first implication is trivial. The other one follows from (i) and the fact that  the predicate ``${\rm cf}(\alpha)=\omega$'' is  $\Sigma_1$.

(iii) Immediate from (ii). \telos

\begin{Cor} \label{C:immediate}
${\rm LZFC}^{\Sigma_1sext}\vdash$ V-P-H. Therefore
$${\rm ZFC}\subseteq {\rm LZFC}^{\Sigma_1sext}+VP\subseteq {\rm LZFC}^{sext}+VP.$$
\end{Cor}

{\em Proof.}  By \ref{L:specsigma} (iii) ${\rm LZFC}^{\Sigma_1sext}$ is stronger than ${\rm LZFC}^{prec}$. So the claim follows from  \ref{T:special}.  [Alternatively, the proof  follows from \ref{L:specsigma} (i). Observe that  the statement ``$(A,R)$ is rigid''  itself is $\Pi_1$. So if $(A,R)\in M$, $M$ is $\Sigma_1$-strongly extendible and $M\models$``$(A,R)$ is rigid'', then $(A,R)$ is absolutely rigid.] \telos

\begin{Cor} \label{Sigma1}
(i) ${\rm LZFC}^{\Sigma_1sext}\vdash$ {\rm Powerset}.

(ii) In ${\rm LZFC}^{\Sigma_1sext}$,  we have $V=\bigcup_\alpha V_\alpha$, where $V_\alpha={\cal P}^\alpha(\emptyset)$. If $M$ is $\Sigma_1$-strongly extendible, $M=V_\alpha$ for some ordinal $\alpha$.
\end{Cor}

{\em Proof.} As in the proofs of \ref{T:surprise} (i) and \ref{P:rankinside}, both claims  just need the fact that every set belongs to a model $M$ such that $M\prec_{\Sigma_1} V$, which holds  in  ${\rm LZFC}^{\Sigma_1sext}$ in view of \ref{L:specsigma} (i). \telos

\vskip 0.2in

We see that the theories ${\rm LZFC}^{spec}$, ${\rm LZFC}^{\Sigma_1sext}$ and ${\rm LZFC}^{sext}$ are reasonable examples of the theory $T$ of \ref{T:main} which, when augmented with $VP$,   restores ZFC. Moreover, as shown in Proposition \ref{P:nostrong}, the consistency strength of   ${\rm LZFC}^{sext}$ (and thus of the rest of the preceding theories) is no higher than that of the existence of a strongly inaccessible cardinal.

\vskip 0.2in

We can prove further that for $T$ being some of the preceding  theories $T+VP$ and ${\rm ZFC}+VP$ are identical. Actually  it suffices to show that ${\rm ZFC}+VP\vdash Loc^{sext}({\rm ZFC})$. For this we need to employ {\em extendible} cardinals the definition of which we recall
(see \cite[p. 379f]{Je03}).

\begin{Def} \label{D:cardextend}
{\em A cardinal $\kappa$ is} extendible {\em if for every $\alpha>\kappa$ there is a $\beta$ and an elementary embedding $j:V_\alpha\rightarrow V_\beta$ such that $\kappa={\rm crit}(j)$. }
\end{Def}

Extendible cardinals are very large. Every extendible cardinal is supercompact, therefore it is strong. So by Lemma  \ref{L:strong},  if $\kappa$ is extendible then $V_\kappa$ is strongly extendible.

\begin{Lem} \label{L:jech}
 In ${\rm ZFC}+VP$ there is a proper class of extendible cardinals.
\end{Lem}

{\em Proof.} See Lemma 20.25 of \cite{Je03} and the remark immediately after its  proof. \telos

\begin{Thm} \label{T:extcard}
${\rm ZFC}+VP\vdash Loc^{sext}({\rm ZFC})$. Hence $${\rm LZFC}^{sext}+VP\subseteq {\rm ZFC}+VP.$$
\end{Thm}

{\em Proof.}  Assume ${\rm ZFC}+VP$ and let $x$ be a set. By \ref{L:jech} there is an extendible cardinal $\kappa$ such that $x\in V_\kappa$. Every extendible cardinal is strong; therefore, by Lemma \ref{L:strong},  $V_\kappa$ is a strongly extendible model. Thus every $x$ belongs to a strongly extendible model and so $Loc^{sext}({\rm ZFC})$ is true. \telos

\begin{Thm}
The theories
(a) ${\rm ZFC}+VP$, (b) ${\rm LZFC}^{sext}+VP$, (c) ${\rm LZFC}^{\Sigma_1sext}+VP$, (d) ${\rm LZFC}^{spec}+VP$ are identical.
\end{Thm}

{\em Proof.} By  \ref{T:special},  \ref{L:specsigma} (iii) and \ref{T:extcard} we have
$${\rm ZFC}+VP\subseteq {\rm LZFC}^{spec}+VP\subseteq {\rm LZFC}^{\Sigma_1sext}+VP\subseteq$$$$ {\rm LZFC}^{sext}+VP\subseteq {\rm ZFC}+VP.$$
   \telos

\vskip 0.2in

It follows from \ref{L:jech} and \ref{T:extcard} that $VP$ has  strong consequences  when combined with ZFC. But what when  combined with LZFC? In particular the following are open.

\begin{Quest} \label{Q:import1}
Does ${\rm LZFC}+VP$ prove V-P-H?
\end{Quest}

\begin{Quest} \label{Q:import2}
Does ${\rm LZFC}+VP$ prove $Loc^{sext}({\rm ZFC})$ (or  some of the weaker principles $Loc^{\Sigma_1sext}({\rm ZFC})$, $Loc^{spec}({\rm ZFC})$)?
\end{Quest}

If the answer to either of the preceding questions is affirmative, then LZFC+$VP$ also restores ZFC, so the  results of this section can be  considerably improved.   However what we were able to prove is only the following.

\begin{Prop} \label{P:addition}
${\rm LZFC}+VP\vdash Loc^{ext}({\rm ZFC})$.
\end{Prop}

{\em Proof.} Let $a$ be a set in the universe of LZFC+$VP$. We have to show that there is an extendible  model $M$ of ZFC such that $a\in M$. Let ${\cal U}=\{(M,\in,\{x\}_{x\in TC(\{a\})}):M\in {\cal M}(a)\}$. In view of \ref{L:properclass} above it is easy to see that  ${\cal U}$ is a proper class.  By $VP$ there are $M\neq N$ in ${\cal M}(a)$ and an elementary embedding $$j:(M,\in,\{x\}_{x\in TC(\{a\})})\rightarrow (N,\in,\{x\}_{x\in TC(\{a\})}).$$ If $j=id$, then $M\prec N$; thus  $M$ is extendible. Let $j\neq id$. By the definition of $j$, $j\restr TC(\{a\})=id$. Let $\kappa={\rm crit}(j)$. By Lemma \ref{L:inac} (iii), $M_\kappa$ is the greatest transitive subset of $M$ fixed by $j$, so $TC(\{a\})\subseteq M_\kappa$. Thus  $a\in M_\kappa$. Also, by \ref{L:inac} (v), $M_\kappa\models{\rm ZFC}$ and, by  \ref{L:inac} (vi), $M_\kappa\prec j(M_\kappa)=N_{j(\kappa)}$. Thus $M_\kappa$ is an extendible model and $a\in M_\kappa$. \telos

\vskip 0.2in

{\bf Acknowledgement} I would like to thank  the anonymous referees for catching several errors in previous drafts  of the paper.


\begin{thebibliography}{99}
\bibitem{AR94}
    J. Ad\'{a}mek and J. Rosick\'{y}, {\em Locally representable and accessible categories,} London Mathematical Society Lecture Note Series, Vol. 189, Cambridge University Press, 1994.
\bibitem{Co06}
   P. Corazza, The spectrum of elementary embeddings $j:V\rightarrow V$, {\em Ann. Pure and Appl. Logic} {\bf 139} (2006), 327-399.
\bibitem{De84}
    K. Devlin, {\em Constructibility,} Perspectives in Mathematical Logic,  volume 6, Springer 1984.
\bibitem{Je03}
    T. Jech, {\em Set Theory,}  the Third Millennium Edition, Springer 2003.
\bibitem{Ka97}
    A. Kanamori, {\em The Higher Infinite,} Perspectives in Mathematical
    Logic, Springer 1997.
\bibitem{MV59}
    R. Montague and R.L. Vaught, Natural models of set theories,
    {\em Fund. Math.} {\bf XLVII} (1959), 219-242.
\bibitem{Tz10}
    A. Tzouvaras, Localizing the axioms, {\em Arch. Math.  Logic} {\bf 49} (2010), no. 5, 571-601, and Erratum to ``Localizing the axioms'', {\em Arch. Math.  Logic} {\bf 50} (2011), no. 3, 513.
\bibitem{VPH65}
    P. Vop\v{e}nka, A. Pultr and Z. Hedrl\'{i}n, A rigid relation exists on any set, {\em Comment. Math. Univ. Carolinae} {\bf 6} (1965), 149-155.

\end{thebibliography}
\end{document}